\documentclass[12pt, a4paper]{article}
\usepackage{setspace}
\usepackage{geometry}
\usepackage[OT1]{fontenc}
\usepackage{amsthm,amsmath}
\usepackage[colorlinks,citecolor=blue,urlcolor=blue]{hyperref}
\usepackage{amssymb,amsmath,bm,mathrsfs,makeidx,amsfonts,graphicx,amsthm,multirow}
\usepackage{booktabs}

\newcommand\X{\mathbf X}

\newcommand\x{\mathbf x}

\renewcommand\Im{{\rm\,Im} }

\newcommand\bR{{\mathbb R}}
\newcommand\Z{{\mathbf Z}}
\newcommand\z{{\mathbf z}}
\newcommand\y{{\mathbf y}}

\newcommand\wh{\widehat }
\newcommand\wt{\widetilde }

\newcommand\bN{{\mathbb N}}
\newcommand\pto{\stackrel p\rightarrow}

\newcommand\e{{\mathbb E}}
\newcommand\p{{\mathbb P}}
\newcommand\var{{\rm Var}}
\newcommand\tr{{\mathrm{tr}}}
\newcommand\cF{{\mathcal F}}
\newcommand\cG{{\mathcal G}}

\newtheorem{theorem}{Theorem}[section]%
\newtheorem{lemma}[theorem]{Lemma}%
\newtheorem{proposition}[theorem]{Proposition}%
\newtheorem{remark}[theorem]{Remark}%

\begin{document}

\begin{center}
\Large The necessary and sufficient conditions in\\ the Marchenko-Pastur theorem.  
\end{center}
\begin{center}
\large Pavel~Yaskov\footnote{National University of Science and Technology MISIS, Russia
\\
Steklov Mathematical Institute, Russia\\
 e-mail: yaskov@mi.ras.ru\\Supported
    by RNF grant 14-21-00162 from the Russian Science Foundation.}
 \end{center}

\begin{abstract} We show that a weak concentration property for quadratic forms of isotropic random vectors $\x$ is necessary and sufficient for the validity of the Marchenko-Pastur theorem for sample covariance matrices  of random vectors having the form $C\x$, where $C$ is any rectangular matrix with orthonormal rows. We also obtain some general conditions  guaranteeing the weak concentration property.
\end{abstract}

\begin{center}
{\bf Keywords:} random matrices; covariance matrices; the Marchenko-Pastur theorem. 
\end{center}

\section{Introduction}
The  Marchenko-Pastur (MP) theorem  \cite{MP} is one of the key results in the random matrix theory. It states that, with probability one, the empirical spectral distribution of a {\it sample covariance matrix}  
\begin{equation}
\label{e8}
\wh \Sigma_n=\frac{1}{n}\sum_{k=1}^n \x_{pk}\x_{pk}^\top
\end{equation}
 weakly converges to the MP law   with parameter $\rho>0$  as $n\to\infty$ if  $p=p(n)$ satisfies $p/n\to \rho$ and, for each $p,$ $\{\x_{pk}\}_{k=1}^n$ are i.i.d. copies of an isotropic $\bR^p$-valued  random vector $\x_p$ satisfying certain assumptions (for definitions, see Section 2).

In the classical case, entries of $\x_p$   are assumed to be i.i.d. copies of some random variable  with mean zero and unit variance (e.g., see Theorem 3.6 in \cite{BS}). In general, entries  of $\x_p=(X_{p1},\ldots,X_{pp}) $ can be any independent random variables that have mean zero and unit variance and satisfy  the Lindeberg condition 
\begin{equation}\label{lin}
L_p(\varepsilon)=\frac{1}{p}\sum_{k=1}^p\e X_{pk}^2 I(|X_{pk}|>\varepsilon \sqrt{p})\to 0\quad\text{for all }\varepsilon>0\text{ as }p\to\infty
\end{equation}
(see \cite{P}). The independence assumption can be relaxed in a number of ways. E.g., in \cite{PP},  the MP theorem is proven for isotropic random vectors $\x_p$ with a centred log-concave distribution.   

All mentioned assumptions imply that quadratic forms $\x_p^\top A_p\x_p$ concentrate near their expectations up to an error term $o( p)$ with probability $1-o(1)$  when $A_p$ is any $p\times p$ real matrix with the spectral norm $\|A_p\|=O(1)$. In fact, this condition is sufficient for the Marchenko-Pastur theorem  (see \cite{BZ}, \cite{G}, \cite{PP}, Theorem 19.1.8 in \cite{PS}, and \cite{Y}).
Recently, it was also proved in \cite{CT} that extreme eigenvalues of $\wh\Sigma_n$ converge in probability to the edges of the  support of the limiting Marchenko-Pastur law if a form of the concentration property for quadratic forms holds (for details, see \cite{CT}).

However, as noted in \cite{A}, the above concentration property for quadratic forms   is {\it not} necessary in general. Namely, take $p=2q$ for  $q=q(n)$ and consider \[\x_p=\sqrt{2}(\z_q \xi, \z_q (1-\xi)),\]
where $\z_q$ is a standard normal vector in $\bR^q$, $\xi$ is a random variable independent of $\z_q$,  and $\p(\xi=0)=\p(\xi=1)=1/2$. Then $\wh\Sigma_n$ is a 2-block-diagonal matrix such that each block satisfies the MP theorem. It can be directly checked that $(\wh\Sigma_n)_{n=1}^\infty$ also satisfy the MP theorem and  each $\x_p$ is an isotropic random vector for which the above concentration property doesn't hold.    

There are many results related to the MP theorem, where some other dependence assumptions are considered. E.g., see  \cite{A}, \cite{BM}, \cite{mp}, \cite{GNT},  \cite{PM}, \cite{MPP}, \cite{OR}, \cite{PSc}, \cite{Y}.

In the present paper we consider the case when not only $\x_p$  but also $C_q\x_p$ satisfies  the MP theorem for each sequence $(C_q)_{p=1}^\infty$, where $q=q(p)\leqslant p$  and  $C_q$ is a $q\times p$ matrix with $C_qC_q^\top =I_q$ for the $q\times q$ identity matrix $I_q$. We prove that the weak concentration property for quadratic forms is a necessary and sufficient assumption in this case. In addition, we derive this property under quite general assumptions recently studied in \cite{PM}. We also show that this property implies some other results in the random matrix theory beyond the MP theorem.

The paper is structured as follows. Section 2 contains some preliminaries, main assumptions, and notation. Main results are presented in Section 3. Section 4 deals with proofs. Some additional results are given in Appendices.

\section{Preliminaries and notation}
We now introduce assumptions and notation that will be used throughout the paper.

 For each $p\geqslant 1,$ let $\x_p$ be a random vector in $\bR^p$.  We call $\x_p$  isotropic if $\e\x_p\x_p^\top=I_p$ for the $p\times p$ identity matrix $I_p$. Let also $ \X_{pn}$ be a $p\times n$ matrix with columns $\{\x_{pk}\}_{k=1}^n$ that are i.i.d. copies of $\x_p$, unless otherwise stated. Then, for $\Sigma_n$ given in \eqref{e8},
\[\wh\Sigma_n=\frac{1}{n}\X_{pn}\X_{pn}^\top.\]

Define also the MP law $\mu_\rho$ with parameter $\rho>0$ by 
\[d\mu_\rho=\max\{1-1/\rho,0\}\,d\delta_0+\frac{\sqrt{(b-x)(x-a)}}{2\pi x\rho   }I(x\in [a,b])\,dx, \]
where $\delta_c$ is a Dirac function with mass at $c$, $a=(1-\sqrt{\rho})^2,$ and $b=(1+\sqrt{\rho})^2.$ Set also  $\mathbb C^+=\{z\in\mathbb C:\,\Im(z)>0\}$.

For a real symmetric $p\times p$ matrix  $A$  with eigenvalues $\lambda_1,\ldots,\lambda_p$, its empirical spectral distribution  is defined by \[\mu_A=\frac{1}{p}\sum_{k=1}^p\delta_{\lambda_k}.\]
If, in addition, $A$ is  positive semi-definite, then $A^{1/2}$ will be the principal square root of $A$. If $A$  is a complex  rectangular  matrix, then $\|A\|$ will be the spectral norm of $A$, i.e. $\|A\|=(\lambda_{\max}(A^*A))^{1/2}$, where  $A^*=\overline{A}^\top$  and $\lambda_{\max}$ denotes the maximal eigenvalue. All matrices below will be real, unless otherwise stated. Let also  $\|v\|$ be the Euclidean norm of   $v=(v_1,\ldots,v_p)\in\mathbb C^p$, i.e. $\|v\|=\big(\sum_{i=1}^p|v_i|^2\big)^{1/2}.$ 

Consider the following assumptions.

$(\mathrm A0)$ $p=p(n)$ satisfies $p/n\to \rho$ for some $\rho>0$ as $n\to\infty$.

$(\mathrm A1)$  $(\x_p^\top A_p\x_p-\tr( A_p))/p\pto 0$ as $p\to\infty$ for all sequences of real symmetric positive semi-definite $p\times p$  matrices $A_p$ with uniformly bounded spectral norms $\|A_p\|$. 

Assumption $(\mathrm A1)$ is a form of the weak law of large numbers for quadratic forms. Stronger forms of $(\mathrm A1)$   (with convergence in $L_2$ instead of convergence in probability) studied in the papers \cite{BZ}, \cite{PP}, and in the book of \cite{PS} (see Chapter 19). In the special case of isotropic $\x_p$,   $\e (\x_p^\top A_p\x_p)=\tr( A_p)$ and $(\mathrm A1)$  states  that $\x_p^\top A_p\x_p$ concentrates near its expectation up to a term $o(p)$ with probability $1-o(1)$ when $\|A_p\|=O(1)$. It  is proved in \cite{Y1} (see also Lemma \ref{l2} in Section 4) that 
\begin{align} \label{equiv}
&\text{if each $\x_p$ has independent entries with mean zero and unit variance,}\nonumber\\
&\text{then $(\mathrm A1)$ is equivalent to the Lindeberg condition \eqref{lin}.}
\end{align}  For general isotropic $\x_p$, we can equivalently reformulate $(\mathrm A1)$ as follows (for a proof, see Appendix A). 
\begin{proposition}\label{p2}
If $\x_p,$ $p\geqslant 1,$ are isotropic, then $(\mathrm A1 )$ holds iff $(\mathrm A1^*)$ holds, where 
 \begin{quote}
$(\mathrm A1^*)$  $(\x_p^\top \Pi_p\x_p-\tr( \Pi_p))/p\pto 0$ as $p\to\infty$ for all sequences of  $p\times p$ orthogonal projection matrices matrices $\Pi_p$. 
 \end{quote}
\end{proposition}

We  will also need a more general form of  $(\text{A1})$ designed for non-isotropic $\x_p$ (namely, when $\e( \x_p\x_p^\top)=\Sigma_p\neq I_p$). For each $p\geqslant 1,$ let $\Sigma_p$ be  a $p\times p$ symmetric positive semi-definite matrix $\Sigma_p$. Consider the following assumptions.

$(\text{A2})$  $(\x_p^\top A_p\x_p-\tr(\Sigma_p A_p))/p\pto 0$ as $p\to\infty$ for all sequences of real symmetric positive semi-definite $p\times p$  matrices $A_p$ with uniformly bounded spectral norms $\|A_p\|$. 

$(\text{A3})$  $\tr(\Sigma_p^2)/p^2\to 0$ as $p\to\infty$. 
\\
The next proposition shows that  $(\text{A2})$ and $(\text{A3})$ are equivalent in the Gaussian case (for a proof, see Appendix A). 
\begin{proposition}\label{p1}
For each $p\geqslant 1,$ let $\x_p$ be a Gaussian vector with mean zero and variance $\Sigma_p$. Then $(\mathrm A2)$ holds if and only if $(\mathrm A3)$ holds. 
\end{proposition} 

We now give a particular and quite general example of $\x_p,$ $p\geqslant 1,$ satisfying $(\mathrm A2)$.
\begin{proposition}\label{p3}
For each $p\geqslant 1,$ let $\x_p=(X_{p1},\ldots,X_{pp})$ be a random vector with mean zero and variance $\Sigma_p$. Suppose there is a nonincreasing sequence $\{\Gamma_j\}_{j=0}^\infty$ such that  $\Gamma_j\to0, $  $j\to\infty$, 
\[\e|\e(X_{pk}|\cF_{k-j}^p)|^2\leqslant \Gamma_j\quad \text{and}\quad \e|\e(X_{pk}X_{pl}|\cF_{k-j}^p)-\e(X_{pk}X_{pl})|\leqslant \Gamma_j\]
for all $1\leqslant k\leqslant l\leqslant p$ and $j=0,\ldots,k$, where $\cF_l^p=\sigma(X_{pk},k\leqslant l),$ $l\geqslant 1,$  and $\cF_0^p$ is the trivial $\sigma$-algebra. If \eqref{lin} holds, then $(\mathrm A2)$ holds. 
\end{proposition} 

 Proposition \ref{p3} is closely related to Theorem 5 in \cite{PM}, where the same dependence conditions are considered, but another result is proven. We will prove Proposition \ref{p3} using Lindeberg's method and Bernstein's block technique as in the proof of Theorem 5 in \cite{PM} (see Appendix A).

Let us also give other versions of $(\text{A2})$ and $(\text{A3})$  allowing some dependence and heterogeneity in $\x_{pk}$ over $k$.  For $\bR^p$-valued random vectors $\{\x_{pk}\}_{k=1}^n$, let $\cF_0^p$ be the trivial $\sigma$-algebra and $\cF_{k}^p=\sigma(\x_{pl},1\leqslant l\leqslant k)$. For given symmetric positive  semi-definite $p\times p$ matrices $\{\Sigma_{pk}\}_{k=1}^n$, introduce the following assumption.

$(\mathrm A2^*)$  For all $\varepsilon>0$ and every stochastic array $\{A_{pk},p\geqslant k\geqslant  1\}$ with symmetric positive semi-definite symmetric $\cF_{k-1}^p$-measurable random $p\times p$ matrices  $A_{pk}$ having $\|A_{pk}\|\leqslant 1$ a.s.,
\[\frac{1}{n}\sum_{k=1}^n\p\big(|\x_{pk}^\top A_{pk}\x_{pk}-\tr(\Sigma_{pk}A_{pk})|>\varepsilon p\big)\to0.\]

$(\mathrm A3^*)$   $(np^2)^{-1}\sum_{k=1}^n\tr(\Sigma_{pk}^2)\to 0$ as $p,n\to\infty$.

In fact, $(\mathrm A2^*)$ and $(\mathrm A3^*)$ are  {\it average} versions of $(\mathrm A2)$ and $(\mathrm A3)$, respectively. In particular, if  $\x_{pk}$ has independent centred entries and variance $\Sigma_{pk}$ for all $p,k$,  then $(\mathrm A2^*)$ follows from   $(\mathrm A3^*)$ and the (second) Lindeberg condition
\begin{equation*}\label{lin2}
\frac{1}{np}\sum_{k=1}^n\sum_{j=1}^p\e X_{kj}^2 I(|X_{kj}|>\varepsilon \sqrt{p})\to 0\quad\text{for all\quad}\varepsilon>0
\end{equation*} 
as $p,n\to\infty$, where $X_{kj}=X_{kj}(p),$ $j=1,\ldots,p$, are entries of $\x_{pk}$. The latter can be checked directly as  in the proof of Proposition 2.1 in \cite{Y1}.

Finally, we introduce the key limit  property for a sequence of random matrices. Let $p=p(n)$ be such that $p/n\to\rho>0 $   and, for each $n,$ let $M_n$ be a symmetric positive semi-definite $p\times p$  matrix. We say that $(M_n)_{n=1}^\infty$  satisfies $(\mathrm {MP})$ if,
\begin{quote} with probability one,  $\mu_{C_qM_nC_q^\top }$ weakly converges to $\mu_{\rho_1}$
as $n\to\infty$ for all $\rho_1\in(0,\rho]$ when $q=q(n)\leqslant p(n)$ satisfies $q/n\to \rho_1$ and $(C_q)_{n=1}^\infty$ is any sequence of matrices such that the size of $C_q$ is $q\times p$ and $C_q C_q^\top=I_q$.
\end{quote}

\section{Main results}

First, we derive necessary and sufficient conditions in the classical setting. 

\begin{theorem}\label{t3}  For each $p\geqslant 1,$ let  $\x_p$ be  isotropic and have centred independent entries. If $p=p(n)$ satisfies $(\mathrm A0)$,  then  $ \mu_{\wh\Sigma_n} $ weakly converges to $\mu_\rho$ almost surely as $n\to\infty$ iff \eqref{lin} holds for given $p=p(n)$. 
\end{theorem}

This result is not new. As far as we know, it was initially established by Girko via a different and less transparent method (e.g., see Theorem 4.1 in \cite{G1}). In our proof of Theorem \ref{t3}, the necessity part  follows from the lemma below.

\begin{lemma}\label{ml1}
Let $(\mathrm A0)$ hold and $\x_p$ be isotropic for all $p=p(n)$. If, with probability one, $\mu_{\wh \Sigma_n}$ weakly converges to $\mu_\rho$ as $n\to\infty$, then, for given $p=p(n)$, 
\begin{equation}\label{e9}
  \frac{ \x_p^\top \x_p-1}{p}\pto 0,\quad n\to\infty.
\end{equation}
\end{lemma}
  
The classical independence setting differs a lot from the general case of isotropic distributions. Namely, when entries of each $\x_p$ are independent and orthonormal,  \eqref{e9} is equivalent to $(\mathrm A1)$ (see the proof of Theorem \ref{t3}). In general,   \eqref{e9}  doesn't imply   $(\mathrm A1)$ as in the counterexample from the Introduction. To get   $(\mathrm A1)$, we need more  than convergence of $\mu_{\wh\Sigma_n}$, e.g., $(\mathrm{MP})$.

We now state the main result of this paper.

\begin{theorem}\label{t2}  Let $(\mathrm A0)$ hold and $\wh\Sigma_n,$ $n\geqslant 1,$ be as in \eqref{e8}. If $(\mathrm A1)$ holds, then $(\wh\Sigma_n)_{n=1}^\infty$ satisfies $(\mathrm{MP})$. Conversely, if  the latter holds and $\x_p$ is isotropic for all $p=p(n)$, then $(\mathrm A1)$ holds for  given $p=p(n)$. 
\end{theorem}

The proof of the necessity part in Theorem \ref{t2} follows from Lemma \ref{ml1}. The sufficiency part can be proved directly as in \cite{Y}, but we prefer to derive it from a more general result.  

\begin{lemma}\label{ml}
Let $(\mathrm A2)$ and $(\mathrm A3)$ hold for some $\Sigma_p,$ $p\geqslant 1$. If $(A0)$ holds for some $p=p(n)$, then, for all $z\in\mathbb C^+$,
\begin{equation}\label{st}
\lim_{n\to\infty}\frac{1}{p}\Big[\tr\big(n^{-1}\wh\X_{pn}\wh\X_{pn}^\top+B_p-zI_p\big)^{-1}- \tr\big(n^{-1}\wh\Z_{pn}\wh\Z_{pn}^\top+B_p-zI_p\big)^{-1}\Big]=0
\end{equation} with probability one, where, for each $n\geqslant 1,$ \begin{itemize}
\item[] $\wh\X_{pn}=\X_{pn}+C_{pn}$ and $ \wh\Z_{pn}=\Z_{pn}+C_{pn},$
\item[] $\Z_{pn}$ is a $p\times n$  matrix with i.i.d. centred Gaussian columns having variance $\Sigma_p$,
\item[] $B_p$ is a $p\times p$  nonrandom   matrix with $\|B_p\|=O(1),$
\item[] $C_{pn}$ is a $p\times n$  nonrandom   matrix such that $\|n^{-1}C_{pn}C_{pn}^\top\|=O(1).$ 
\end{itemize}
\end{lemma}
  
\begin{remark}\label{r2}{\normalfont  By the definition of $\mu_{A}$,
\[\int_{-\infty}^\infty\frac{ \mu_{A}(d\lambda)}{\lambda -z}=\frac{1}{p} \tr\big(A-zI_p)^{-1},\quad z\in\mathbb C^+,\] where $A$ is a symmetric $p\times p$ matrix. Therefore, by the Stieltjes continuity theorem (e.g., see Exercise 2.4.10 in \cite{Tao}), \eqref{st} implies that if $\mu_{\wh\Sigma_n}$ with $\wh\Sigma_n=n^{-1}\wh\Z_{pn}\wh\Z_{pn}^\top+B_p$ weakly converges to some measure $\mu $ a.s., then $\mu_{\Sigma_n}$ with $\Sigma_n=n^{-1}\wh\X_{pn}\wh\X_{pn}^\top+B_p$ weakly converges to the same measure $\mu $ a.s.}
\end{remark}

We now consider a generalization of Lemma \ref{ml} allowing some dependence and heterogeneity in $\x_{pk}$ over $k$. Let $\X_{pn}$ be a $p\times n$ random matrix with columns $\{\x_{pk}\}_{k=1}^n$.  Let also  $\{\Sigma_{pk}\}_{k=1}^n$ be symmetric positive semi-definite  $p\times p$ matrices and 
$\Z_{pn}$ be a $p\times n$ Gaussian random matrix whose $k$-th column $\z_{pk}$, $k=1,\ldots,n$, has mean zero and variance $\Sigma_{pk}$. Assume also that $\Z_{pn}$ and $\X_{pn}$ are independent.
\begin{lemma}\label{ml2}
Let $(\mathrm A2^*)$, $(\mathrm A3^*)$, and $(\mathrm A0)$ hold for some $p=p(n)$. Then, for all $z\in\mathbb C^+$, \eqref{st} with convergence in probability holds,  where  $B_p$, $C_{pn}$, $\wh\X_{pn}$, and $\wh \Z_{pn}$ are  as in Lemma \ref{ml}. 
\end{lemma}

\section{Proofs} 
\noindent{\bf Proof of Theorem \ref{t3}.} By \eqref{equiv}, the Lindeberg method \eqref{lin} is equivalent to $(\mathrm A1)$ (see also Lemma \ref{l2} below). However, the proof of  this proposition is still valid without the assumption that $\x_p$ has mean zero. Now, if $(\mathrm A1)$ holds, then $\mu_{\wh\Sigma_n}$ weakly converges to $\mu_\rho$ almost surely as $n\to\infty$  by Theorem \ref{t2}. 

Suppose $\mu_{\wh\Sigma_n}$ weakly converges to $\mu_\rho$ a.s. We need to prove \eqref{lin}. By Lemma \ref{ml1}, \eqref{e9} holds.  By the Gnedenko-Kolmogorov necessary and sufficient conditions for relative stability (see (A) and (B) in \cite{H}), \eqref{lin} is equivalent to \eqref{e9}. This finishes the proof of the theorem. Q.e.d.

\noindent{\bf Proof of Lemma \ref{ml1}.}  We will proceed as in \cite{Y1}. Let $\x_p=\x_{p,n+1}$ be independent of the matrix $\X_{pn}$ and distributed as its columns $\{\x_{pk}\}_{k=1}^n$. Define also \[A_n=n\wh\Sigma_n=\X_{pn}\X_{pn}^\top=\sum_{k=1}^n\x_{pk}\x_{pk}^\top\quad \text{and }\quad B_n=A_n+\x_{p}\x_{p}^\top=\sum_{k=1}^{n+1}\x_{pk}\x_{pk}^\top.\]
Fix $\varepsilon>0.$ The matrix $B_n +\varepsilon nI_p$ is non-degenerate and
\[p=\tr\big((B_n +\varepsilon nI_p)(B_n +\varepsilon nI_p)^{-1}\big)=\sum_{k=1}^{n+1}\x_{pk}^\top(B_n +\varepsilon nI_p)^{-1}\x_{pk}+\varepsilon n\,\tr(B_n +\varepsilon nI_p)^{-1}.\]
Taking expectations and using the exchangeability of $\{\x_{pk}\}_{k=1}^{n+1}$,
\begin{align}\label{e1}
p=&(n+1)\e \x_p^\top(B_n +\varepsilon nI_p)^{-1}\x_p+\varepsilon n\,\e \tr(B_n +\varepsilon nI_p)^{-1}.
\end{align}

Define $f_n(\varepsilon)=\tr(A_n+\varepsilon nI_p)^{-1}$. By Facts 5 and 7 in Appendix C,
\[\e\tr(B_n+\varepsilon nI_p)^{-1}=\e f_n(\varepsilon)+o(1)\quad\text{and}\quad \e \x_p^\top(B_n +\varepsilon nI_p)^{-1}\x_p=O(1).\]
Thus,
\begin{equation}
\label{e3}
p/n=\e \x_p^\top(B_n +\varepsilon nI_p)^{-1}\x_p+\varepsilon \e \tr(A_n +\varepsilon nI_p)^{-1}+o(1).
\end{equation}

Let now  $\Z_{pn}$ be a $p\times n$ matrix with i.i.d. $\mathcal N(0,1)$ entries. By the MP theorem,  $\mu_{\widetilde\Sigma_n}$ weakly converges to $\mu_\rho$ a.s., where  $\widetilde\Sigma_n=n^{-1}C_n$ and $C_n=\Z_{pn}\Z_{pn}^\top$. Therefore, by the Stieltjes continuity theorem (see Theorem B.9 in \cite{BS}),
\[\p(S_n(z)-s_n(z)\to 0 \text{\;\;for all }z\in\mathbb C^+)=1,\]
where $S_n(z)$ and $s_n(z)$ are  the Stieltjes transforms defined by
\[S_n(z)=\int_{-\infty}^\infty\frac{ \mu_{\wh\Sigma_n}(d\lambda)}{\lambda -z}=p^{-1}\tr\big(\wh\Sigma_n-zI_p\big)^{-1}=p^{-1}\tr\big(n^{-1}A_n-zI_p\big)^{-1},\]
\[s_n(z)=\int_{-\infty}^\infty\frac{ \mu_{\wt\Sigma_n}(d\lambda)}{\lambda -z}=p^{-1}\tr\big(\wt \Sigma_n-zI_p\big)^{-1}=p^{-1}\tr\big(n^{-1}C_n-zI_p\big)^{-1}.\] 
By Fact 8 in Appendix C and $p/n\to \rho >0,$ the latter implies that 
\begin{equation}
\label{e101}
\p(\tr(A_n +\varepsilon nI_p)^{-1}-\tr(C_n +\varepsilon nI_p)^{-1}\to 0\text{\;\;for all }\varepsilon>0)=1.\end{equation}
Therefore, by the dominated convergence theorem, 
\[
\e \tr(A_n +\varepsilon nI_p)^{-1}=\e \tr(C_n +\varepsilon nI_p)^{-1}+o(1)\text{\;\;for all }\varepsilon>0.
\]

In addition, arguing as above, we derive
\begin{equation}
\label{e4}
p/n=\e \z_p^\top(D_n +\varepsilon nI_p)^{-1}\z_p+\varepsilon \e \tr(C_n +\varepsilon nI_p)^{-1}+o(1),
\end{equation}
where $D_n=C_n+\z_p\z_p^\top$ and $\z_p$ is a $\bR^p$-valued  standard normal vector independent of $C_n$. Subtracting \eqref{e3} from \eqref{e4}, we get
\[\e \x_p^\top(B_n +\varepsilon nI_p)^{-1}\x_p=\e \z_p^\top(D_n +\varepsilon nI_p)^{-1}\z_p+o(1)\quad\text{for all }\varepsilon>0.\]
By Fact 7 in Appendix C and \eqref{e5} (see Appendix B), 
\[\e \z_p^\top(D_n +\varepsilon nI_p)^{-1}\z_p=\e \frac{Z_n}{1+Z_n}=\e  \frac{\e (Z_n|C_n)}{1+\e (Z_n|C_n) }-R_n,\]
where $Z_n=\z_p^\top(C_n +\varepsilon nI_p)^{-1}\z_p$ and 
\begin{align*}
 R_n=\e\frac{(Z_n-\e (Z_n|C_n))^2}{(1+Z_n)(1+\e (Z_n|C_n))^2}&\leqslant \e (Z_n-\e( Z_n|C_n))^2=\e \var(Z_n|C_n)
 \end{align*}
Since $\z_p$ and $C_n$ are independent, $\e (Z_n|C_n)=\tr(C_n +\varepsilon nI_p)^{-1},$ and \[\|(C_n +\varepsilon nI_p)^{-1}\|\leqslant (\varepsilon n)^{-1},\]
we have $\e \var(Z_n|C_n)=o(1)$ by \eqref{v1} (see Appendix A). Hence, $R_n\to 0.$

By \eqref{e101} and the dominated convergence theorem,  
\begin{align*}
\e\frac{\e(Z_n|C_n)}{1+\e(Z_n|C_n)}&=\e\frac{ \tr(A_n +\varepsilon nI_p)^{-1}}{1+ \tr(A_n +\varepsilon nI_p)^{-1}}+o(1).
\end{align*}
Combining the above relations with Fact 7 in Appendix C yields
\begin{align*}\e \x_p^\top(B_n +\varepsilon nI_p)^{-1}\x_p&=\e\frac{ \x_p^\top(A_n +\varepsilon nI_p)^{-1}\x_p}{1+\x_p^\top(A_n +\varepsilon nI_p)^{-1}\x_p}\\
&=\e\frac{ \tr(A_n +\varepsilon nI_p)^{-1}}{1+ \tr(A_n +\varepsilon nI_p)^{-1}}+o(1).
\end{align*}
Additionally,  $\e(\x_p^\top (A_n +\varepsilon nI_p)^{-1}\x_p |A_n)=\tr(A_n +\varepsilon nI_p)^{-1}$ a.s. by the independence of $\x_p$ and $A_n.$ 
\begin{lemma} \label{l4}
For each $n\geqslant 1,$ let $Z_n$ be a random variable such that $Z_n\geqslant 0$ a.s. and $\e Z_n$ is bounded over $n$. If $Y_n,$ $n\geqslant 1,$ are such random elements that 
\[\e\frac{Z_n}{1+Z_n}-\e\frac{\e (Z_n|Y_n) }{1+\e (Z_n|Y_n)}\to 0,\quad n\to\infty,\]
then $Z_n-\e(Z_n|Y_n) \pto0.$
\end{lemma}
The proof of Lemma \ref{l4} is given in Appendix B. Using Lemma \ref{l4}, we conclude that
$\x_p^\top(A_n +\varepsilon nI_p)^{-1}\x_p-\tr(A_n +\varepsilon nI_p)^{-1}\pto 0.$
Multiplying by $\varepsilon$ and $n/p,$ we finally arrive at 
\[p^{-1}(\x_p^\top(\varepsilon^{-1}\wh\Sigma_n + I_p)^{-1}\x_p-\tr(\varepsilon^{-1}\wh\Sigma_n +I_p)^{-1})\pto 0\quad\text{for all $\varepsilon>0,$}\]
where $\wh\Sigma_n=A_n/n$ is independent of $\x_p$.
Thus, we can find $\varepsilon_n$ that slowly tend to infinity and are such that
\[J_n=p^{-1}(\x_p^\top(\varepsilon_n^{-1}\wh\Sigma_n + I_p)^{-1}\x_p-\tr(\varepsilon_n^{-1}\wh\Sigma_n +I_p)^{-1})\pto 0.\]

We know that  $\mu_{\wh\Sigma_n}$ weakly converges to $\mu_\rho$ a.s. The support of $\mu_\rho$ is bounded. Hence,
writing $\varepsilon_n^{-1}\wh\Sigma_n=\sum_{k=1}^p\lambda_ke_ke_k^\top$ for some $\lambda_k=\lambda_k(n)\geqslant0$ and  orthonormal vectors $e_k=e_k(n)\in\bR^p$, $k=1,\ldots,p,$ we conclude that  
\[\frac{1}{p}\sum_{k=1}^pI(\lambda_k>\delta_n)\pto0\]
when $\delta_n=K\varepsilon_n^{-1}\to0$ and $K>0$ is large enough.  In addition,
\[J_n-\frac{\x_p^\top \x_p-1}{p}=U_n+V_n,\]
where 
\[U_n=\frac{1}{p}\sum_{k:\,\lambda_k\leqslant\delta_n}((\x_p,e_k)^2-1)\bigg(\frac{1}{\lambda_k+1}-1\bigg),\]
\[V_n=\frac{1}{p}\sum_{k:\,\lambda_k>\delta_n}((\x_p,e_k)^2-1)\bigg(\frac{1}{\lambda_k+1}-1\bigg).\]
We finish the proof by showing that $U_n\pto 0$ and $V_n\pto0.$
By the independence of $\wh\Sigma_n$ and $\x_p,$  we have $\e((\x_p,e_k)^2|\wh\Sigma_n)=e_k^\top e_k=1$. In addition, \[\e |U_n|=\e[ \e(|U_n||\wh\Sigma_n)]\leqslant \frac{2}p\,\e\sum_{k:\,\lambda_k\leqslant\delta_n}\frac
{\lambda_k}{\lambda_k+1}\leqslant 2\delta_n=o(1),\]
\[\e |V_n|=\e[ \e(|V_n||\wh\Sigma_n)]\leqslant \frac{2}p\e\sum_{k=1}^p I(\lambda_k>\delta_n)=o(1).\]
Finally, we conclude that $( \x_p^\top \x_p-1)/p=J_n-(U_n+V_n)\pto0.$ Q.e.d.

\noindent{\bf Proof of Theorem \ref{t2}.} Let $(\mathrm A1)$ hold.
Fix some $q=q(n)$ such that $q\leqslant p$ and $q/n\to\rho_1>0.$  For each $q,$ let  $C_q$ be a  $q\times p$ matrix with $C_q C_q^\top=I_q.$  For such $C_q,$   $(\mathrm A1)$ implies that $(\mathrm A2)$ holds for  $(\x_p,\Sigma_p,p)$ replaced by $(C_q\x_p,I_q,q)$.  By Lemma \ref{ml} (without $B_p$ and $C_{pn}$) and Remark \ref{r2}, $\mu_{C_q  \wh\Sigma_n C_q^\top}$ and $\mu_{\widetilde\Sigma_n}$ weakly converge to the same limit  a.s., where $\widetilde\Sigma_n= n^{-1}\Z_{qn}\Z_{qn}^\top$ and $\Z_{qn}$ is a $q\times n$ random matrix with i.i.d. $\mathcal N(0,1)$ entries. By the MP theorem, $\mu_{\widetilde \Sigma_n}$ weakly converges to $\mu_{\rho_1}$  a.s. Thus,  $(\wh\Sigma_n)_{n=1}^\infty$ satisfies $(\mathrm{MP})$.

Let now  $(\wh\Sigma_n)_{n=1}^\infty$ satisfies $(\mathrm {MP})$ and  let $\x_p$ be isotropic for each $p=p(n)$. We need to show that   $(\mathrm A1)$  holds for given $p=p(n)$. By Proposition \ref{p2}, $(\mathrm A1)$ is equivalent to $(\mathrm A1^*)$. Let us verify $(\mathrm A1^*)$.  Fix any  sequence of orthogonal projectors  $\Pi_p, $ $p=p(n),$ such that the size of $\Pi_p$  is $p\times p$. We need to show that
\[\frac{1}{p}\big(\x_p^\top  \Pi_p \x_p-\tr(\Pi_p)\big)\pto0.\]
The latter is equivalent to the fact that 
\[\frac{1}{p_k}\big(\x_{p_k}^\top  \Pi_{p_k} \x_{p_k}-\tr(\Pi_{p_k})\big)\pto0\]
 for any subsequence $\Pi_{p_k}$ such that   $\tr(\Pi_{p_k})/p_k $ has a limit as $k\to\infty$. If this limit is zero, then, obviously, the above convergence  holds.

Assume further w.l.o.g. that  $q/n$ (or, equivalently, $q/p$) has a positive limit, where $q=\tr(\Pi_p)$. Namely, let  $q/n\to \rho_1>0 $ (and $q/p\to\rho_1/\rho$). Write $\Pi_p=C_p^\top D_pC_p ,$ where $C_p$ is a $p\times p$ orthogonal matrix and $D_p$ is a diagonal matrix whose diagonal entries are $1,\ldots,1,0,\ldots,0$ with  $q$ ($=\tr(\Pi_p)$) ones. Let $C_{qp}$ is the $q\times p$ upper block of $C_p$. Then
$\x_p^\top  \Pi_p \x_p=(C_{qp}\x_p)^\top C_{qp}\x_p$.
By $(\mathrm{MP})$, $\mu_{C_{qp}\Sigma_n C_{qp}^\top}$ weakly converges to $\mu_{\rho_1}$ a.s. In addition, $C_{qp}\x_p$ is isotropic, since $\e C_{qp}\x_p(C_{qp}\x_p)^\top=C_{qp}C_{qp}^\top=I_q.$ Hence, by Lemma \ref{ml1}, we get 
\[\frac{1}{q}\big((C_{qp}\x_p)^\top C_{qp}\x_p-q\big)=(\rho/ \rho_1+o(1))\frac{1}{p}\big(\x_p^\top  \Pi_p \x_p-\tr(\Pi_p)\big)\pto0.\]
This proves that $(\mathrm A1^*)$ holds. Q.e.d.

{\bf Proof of Lemma \ref{ml}.} Fix $z\in\mathbb C^+.$ 
Assume w.l.o.g. that $B_p$ is positive semi-definite (we can always replace $B_p$ by $B_p+mI_p$ for $m=\sup_{p}\|B_p\|$ and change $z$ to $z+m$). First we proceed as in Step 1 of the proof of Theorem 1.1 in \cite{BZ} (see also the proof of (4.5.6) on page 83 in \cite{BS}) to show that  $S_n(z)-\e S_n(z)\to 0$ a.s. as $n\to\infty,$  where
\[S_n(z)=p^{-1}\tr\big(n^{-1}\wh\X_{pn}\wh\X_{pn}^\top+B_p-zI_p\big)^{-1}.\]  Similar arguments yield that $s_n(z)-\e s_n(z)\to 0$ a.s. as $n\to\infty,$ where \[s_n(z)=p^{-1}\tr\big(n^{-1}\wh\Z_{pn}\wh\Z_{pn}^\top+B_p-zI_p\big)^{-1}.\] Hence, we only need to show that  $\e S_n(z)-\e s_n(z)\to 0$. 

First, we consider  the case when all entries of $C_{pn}$ are  zeros, i.e. $\wh \X_{pn}=\X_{pn}$.
 Let $\X_{pn}$ and  $\Z_{pn}$ be independent. We will use  Lindeberg's method as in the proof of  Theorem 6.1 in \cite{G}.
 Recall that 
\[\X_{pn}\X_{pn}^\top=\sum_{k=1}^n\x_{pk}\x_{pk}^\top\quad\text{and}\quad
\Z_{pn}\Z_{pn}^\top=\sum_{k=1}^n\z_{pk}\z_{pk}^\top,\]
where $\{\x_{pk}\}_{k=1}^n$ and $\{\z_{pk}\}_{k=1}^n$ are columns of  $\X_{pn}$ and  $\Z_{pn}$, respectively.  If $\z_p$ is a centred Gaussian vector with variance $\Sigma_p$ that is independent of $\x_p,$ then $\{(\x_k,\z_k)\}_{k=1}^n$ are i.i.d. copies of $(\x_p,\z_p)$. In what follows, we  omit the index $p$ and, for example, write $(\x_k,\z_k)$ instead of $(\x_{pk},\z_{pk})$.  

Let us now prove that
\begin{equation}\label{e102}
\e|S_n(z)-s_n(z)|\to 0. 
\end{equation}
Using this representation, we derive that 
\begin{align*}
S_n(z)&=\frac{1}{p}\tr\Big(n^{-1}\sum_{k=1}^n\x_{k}\x_{k}^\top+B_p-z I_p\Big)^{-1},\\
 s_n(z)&=\frac{1}{p}\tr\Big(n^{-1}\sum_{k=1}^n\z_{k}\z_{k}^\top+B_p-z  I_p\Big)^{-1},
\end{align*}
and $|S_n(z)-s_n(z)|\leqslant \sum_{k=1}^n|\Delta_{k}|/p,$ where 
\begin{align*}
\Delta_{k} =\tr\Big(C_{k}+\frac{\x_{k}\x_{k}^\top }n+B_p-z I_p\Big)^{-1}-\tr\Big(C_{k}+ \frac{\z_{k}\z_{k}^\top }{n}+B_p-z I_p\Big)^{-1}
\end{align*}
for
$C_{1}=\sum_{i=2}^n \z_{i}\z_{i}^\top /n,$ $C_{n}=\sum_{i=1}^{n-1} \x_{i}\x_{i}^\top /n,$ and
\begin{equation}\label{e103}
C_{k}=\frac{1}{n}\sum_{i=1}^{k-1} \x_{i}\x_{i}^\top+\frac{1}{n}\sum_{i=k+1}^n \z_{i}\z_{i}^\top ,\quad 1<k<n.
\end{equation}
By Fact 7 in Appendix C,
\[
\tr(C+ww^\top-zI_p)^{-1}-\tr(C-zI_p)^{-1}=-\frac{w^\top (C-zI_p)^{-2}w}{1+w^\top (C-zI_p)^{-1}w}
\] for any real $p\times p$ matrix $C$ and $w\in\bR^p$. Adding and subtracting $\tr(C_{k}+B_p-z I_p)^{-1}$  yield
\[\Delta_{k} =-\frac{ \x_{k}^\top A_{k}^{2}\x_{k}/n}{1+ \x_{k}^\top A_{k}\x_{k}/n}+\frac{ \z_{k}^\top A_{k}^2\z_{k}/n}{1+\z_{k}^\top A_{k}\z_{k}/n},\] where we set $A_{k}=A_{k}(z)=(C_{k}+B_p-zI_p)^{-1}$, $1\leqslant k\leqslant n.$

Let us show that \[\frac1p\sum_{k=1}^n\e|\Delta_{k}|\to 0.\]
The latter implies that $\e |S_n(z)-  s_n(z)|\to 0$.   

Fix some $k\in\{1,\ldots,n\}$. For notational simplicity, we will further write \[\x_p,\z_p,A,\Delta,C\quad \text{instead of}\quad \x_{k},\z_k,A_{k},\Delta_{k},C_{k}\] and use the following properties:   $C+B_p$ is a real symmetric positive semi-definite $p\times p$ random matrix, $(\x_p,\z_p)$ is independent of $A=(C+B_p-zI_p)^{-1}$. 

Fix any $\varepsilon>0$ and let further $v=\Im(z)$ ($ >0$). Take
\[D=\bigcap_{j=1}^2\{|\x_p^\top A^j \x_p-\z_p^\top A^j \z_p|\leqslant \varepsilon p\}\]  
 and derive that $\e |\Delta|\leqslant \e|\Delta|I(D)+2\p\big(\overline{D}\big)/v,$ where    $|\Delta|\leqslant 2/v$ by Fact 5 in Appendix C.

By the law of iterated mathematical expectations and Fact 4 in Appendix C,
\[\p\big(\overline{D}\big)=\e\big[\p\big(\overline{D}|A\big)\big]\leqslant 2\sup_{A_p}\p(|\x_p^\top  A_p \x_p-\z_p^\top A_p \z_p|>\varepsilon  p),\]
where $ A_p$ is any  $p\times p$ complex symmetric matrix  with $\|A_p\|\leqslant M:=\max\{v^{-1},v^{-2}\}$. 

To estimate $\e|\Delta|I(D)$ we need the following technical lemma that is proved in Appendix B.
\begin{lemma} \label{l3} Let $z_1,z_2,w_1,w_2\in\mathbb C$. If $|z_1-z_2|\leqslant \gamma,$ $|w_1-w_2|\leqslant \gamma,$
\[\frac{|z_1|}{|1+w_1|}\leqslant M,\] and $|1+w_2|\geqslant\delta$ for some $\delta,M>0$ and $\gamma\in(0,\delta/2)$, then, 
for some $C=C(\delta,M)>0,$ 
\[\Big|\frac{z_1}{1+w_1}-\frac{z_2}{1+w_2}\Big|\leqslant C\gamma.\]
\end{lemma}

Since $\z_p^\top A\z_p/n=\tr((\z_p\z_p^\top/n) A)$,  Fact 6 in Appendix C implies that
\begin{equation}\label{1w}
|1+\z_p^\top A\z_p/n|\geqslant \delta:=\frac{v}{|z|}.
\end{equation}
Take $\gamma=\varepsilon p/n,$
\[(z_1,w_1)=(\x_p^\top A^2 \x_p, \x_p^\top A \x_p)/n,\qquad 
(z_2,w_2)=(\z_p^\top A^2 \z_p, \z_p^\top A \z_p)/n\]
in Lemma \ref{l3}. By Fact 5 in Appendix C, $ |z_1|/|1+w_1| \leqslant 1/ v .$
By \eqref{1w}, 
\[|1+w_2|\geqslant \delta>\frac{2\varepsilon p}n=2\gamma\]
for small enough $\varepsilon>0$ and large enough $p$ (since $p/n\to \rho>0$).

Using Lemma \ref{l3}, we derive
\[\e|\Delta| I(D)\leqslant C(\delta,1/v)\,\frac{\varepsilon p}{n}.\]
Combining all above estimates together yields
\[ \e |\Delta_{k}| \leqslant C(\delta,1/v)\frac{\varepsilon p}{n}+ \frac 4v\,\sup_{ A_p}\p(|\x_p^\top  A_p \x_p-\z_p^\top  A_p \z_p|>\varepsilon  p)
\]
for each $k=1,\ldots,n$ and
\begin{align}\label{las} \frac1p\sum_{k=1}^n\e |\Delta_{k}| &\leqslant C(\delta,1/v)\varepsilon + \frac {4n}{vp}\,\sup_{ A_p}\p(|\x_p^\top  A_p \x_p-\z_p^\top  A_p \z_p|>\varepsilon  p).\end{align}
Note also that
\begin{align}\label{e106}\p(|\x_p^\top  A_p \x_p-\z_p^\top  A_p \z_p|>\varepsilon  p)\leqslant& \p(|\x_p^\top  A_p \x_p-\tr(\Sigma_pA_p)|>\varepsilon  p/2)+\nonumber\\&\quad+\p(|\z_p^\top  A_p \z_p -\tr(\Sigma_pA_p)|>\varepsilon  p/2).\end{align}
Taking $\varepsilon$ small enough and then $p$ large enough, we can make the right-hand side of \eqref{las}  arbitrarily small by the following lemma (that also holds for $\z_p$ instead of $\x_p$ by $(\mathrm A3)$ and Proposition \ref{p1}).
\begin{lemma} \label{l2} Let $(\mathrm A2)$ holds.  Then, for each $\varepsilon,M>0$,
\begin{equation}\label{limsup}
\lim_{p\to\infty}\sup_{A_p}\p(|\x_p^\top A_p \x_p-\tr(\Sigma_pA_p)|>\varepsilon p)=0,\end{equation}
where the supremum is taken over all complex $p\times p$ matrices $A_p$ with $\|A_p\|\leqslant M.$
\end{lemma}
The proof of Lemma \ref{l2} can be found in Appendix B. Finally, we conclude that 
\[\lim_{n\to\infty}\frac1p\sum_{k=1}^n\e  |\Delta_{k}| =0.\] 
This finishes the proof in the case when all entries of $C_{pn}$ are zeros.

Consider now the case with nonzero $C_{pn}$. Let $ c_k=c_k(n),$ $1\leqslant k\leqslant n,$ be columns of $C_{pn}$. Since $\|n^{-1}C_{pn}C_{pn}^\top\|=\|n^{-1}C_{pn}^\top C_{pn}\|=O(1)$  and $p/n\to\rho>0,$ we have $\max_{1\leqslant k\leqslant n}(c_k^\top c_k)=O(p)$, where $c_k^\top c_k$ are diagonal entries of $n^{-1}C_{pn}^\top C_{pn}$. 

We also have
\[\wh\X_{pn}\wh\X_{pn}^\top=\sum_{k=1}^n\wh\x_{k}\wh\x_{k}^\top\quad\text{and}\quad
\wh\Z_{pn}\wh\Z_{pn}^\top=\sum_{k=1}^n\wh\z_{k}\wh\z_{k}^\top,\]
where $\wh \x_k=\x_k+c_k$ and  $\wh \z_k=\z_k+c_k$ (here $\{(\x_k,\z_k)\}_{k=1}^p$ are i.i.d. copies of $(\x_p,\z_p)$). Arguing  as above,  we conclude that \eqref{las} holds with  
\[\sup_{ A_p}\p(|\x_p^\top  A_p \x_p-\z_p^\top  A_p \z_p|>\varepsilon  p)\] replaced by \[\frac{1}{n}\sum_{k=1}^n\sup_{ A_p}\p(|(\x_p+c_k)^\top  A_p (\x_p+c_k)-(\z_p+c_k)^\top  A_p (\z_p+c_k)|>\varepsilon  p).\]
Recalling that $A_p$ is symmetric, we derive 
\[|(\x_p+c_k)^\top  A_p (\x_p+c_k)-(\z_p+c_k)^\top  A_p (\z_p+c_k)|\leqslant\]\[\leqslant |\x_p^\top  A_p \x_p- \z_p^\top  A_p  \z_p|+2|c_k^\top  A_p  \x_p|+2|c_k^\top  A_p  \z_p|.\]
In addition, 
\[|c_k^\top  A_p  \x_p|^2/p^2\leqslant |\x_p^\top M_k\x_p-\tr(M_k)|/p +|\tr(M_k)|/p \] and the same inequality holds for $\z_p$ instead of $\x_p,$ where $M_k=p^{-1} A_p^* c_k c_k^\top  A_p.$ Note also that, uniformly in $A_p$ with $\|A_p\|\leqslant M,$ 
\[\|M_k\|\leqslant \frac{\|A_p\|^2}p\max_{1\leqslant k\leqslant n}(c_k^\top c_k) =O(1)\] as well as \[\frac{|\tr(M_k)|}p=\frac{| c_k^\top  A_pA_p^* c_k| }{p^2 }\leqslant \frac{\|A_p\|^2}{p^2}\max_{1\leqslant k\leqslant n}  (c_k^\top  c_k) =o(1).\]

Combining these estimates, Lemma \ref{l2},  $(\mathrm A2)$, $(\mathrm A3)$, and Proposition \ref{p1}, we get 
\[\frac{1}{n}\sum_{k=1}^n\sup_{ A_p}\p(|(\x_p+c_k)^\top  A_p (\x_p+c_k)-(\z_p+c_k)^\top  A_p (\z_p+c_k)|>\varepsilon  p) \to 0\]
for all fixed $\varepsilon>0$. Now, we can finish the proof as in the  case with zero $C_{pn}.$   Q.e.d.

\noindent{\bf Proof of Lemma \ref{ml2}.}  As in the proof of Lemma \ref{ml}, we can assume that $B_p$ is positive definite and, first, consider the case with null  $C_{pn}$.  Write further $\x_k$ and $\z_k$ instead of  $\x_{pk}$ and $\z_{pk}$ and note that $(\mathrm A2^*)$ implies that, for all $\varepsilon,M>0,$
\begin{equation}
\label{e105}
\frac{1}{n}\sum_{k=1}^n\e\sup_{A_{pk}}\p(|\x_k^\top  A_{pk} \x_k-\tr(\Sigma_{pk}A_{pk} )|>\varepsilon  p|\cF_{k-1}^p)\to0,
\end{equation}
where the $k$-th supremum is taken over all $\cF_{k-1}^p$-measurable symmetric positive semi-definite random $p\times p$   matrices $A_{pk}$ with $\|A_{pk}\|\leqslant M $ a.s. Recall also that $\cF_{k-1}^p=\sigma(\x_l,l\leqslant k-1)$. Arguing in the same way as in the proof of Lemma \ref{l2}, it can be shown that \eqref{e105} holds even when  the $k$-th supremum is taken over all $\cF_{k-1}^p$-measurable complex random $p\times p$   matrices $A_{pk}$ with $\|A_{pk}\|\leqslant M $ a.s.
Due to $(\mathrm A3^*)$ and \eqref{inequality} (see Appendix A), the same relation holds for $\x_k$ replaced by $\z_k$ (here $\z_k$ is independent of $\cF_{k-1}^p$).

Now, using the same arguments as in the proof of Lemma \ref{ml} after \eqref{e102}, we conclude that \eqref{las} holds with  
\[\sup_{ A_p}\p(|\x_p^\top  A_p \x_p-\z_p^\top  A_p \z_p|>\varepsilon  p)\] replaced by 
\[V_n(\varepsilon):=\frac{1}{n}\sum_{k=1}^n[\p(|\x_k^\top  A_k \x_k-\z_k^\top  A_k \z_k|>\varepsilon  p)+\p(|\x_k^\top  A_k^2 \x_k-\z_k^\top  A_k^2 \z_k|>\varepsilon  p)],\]
where  $A_k=(C_k+B_p-zI_p)^{-1}$ has $\|A_k\|\leqslant M=\max\{1/\Im(z),1/|\Im(z)|^2\}$ a.s.  and  $C_k$ is defined \eqref{e103}.

Each $C_k$ can be written as the sum of two matrices $C_{k1} $ and $C_{k2}$ such that $C_{k1}$ is $\cF_{k-1}^p$-measurable   and $C_{k2}$ 
is a function of $\z_l, l\geqslant k+1$. Since  $\{\z_k\}_{k=1}^n$ are mutually independent and independent from everything else, 
\begin{equation}
\label{e104}
V_n(\varepsilon)\leqslant \frac{2}{n}\sum_{k=1}^n\e\sup_{A_{pk}}\p(|\x_k^\top  A_{pk} \x_k-\z_k^\top  A_{pk} \z_k|>\varepsilon  p|\cF_{k-1}^p),\end{equation}
where the $k$-th supremum is taken over all $\cF_{k-1}^p$-measurable complex $p\times p$ random  matrices $A_{pk}$ having $\|A_{pk}\|\leqslant M $ a.s.  Thus, it follows from \eqref{e106}, \eqref{e105}, and the same relation with $\x_k$ replaced by $\z_k$ that 
\[V_n(\varepsilon)\to0\quad \text{for all}\quad\varepsilon>0.\]
Now, we can finish the proof  as in the proof of Lemma \ref{ml} after \eqref{e106}. Q.e.d.

\section*{Appendix A}
\noindent{\bf Proof of Proposition \ref{p2}.} By definition, $(\mathrm A1^*)$ follows from $(\mathrm A1)$. Suppose $(\mathrm A1^*)$ holds. Let us show that $(\mathrm A1)$ holds. Note that  if  $\Pi_p,$ $p\geqslant1,$ are orthogonal projectors, then  
  $Z_p\pto 1$  and $\e Z_p=1$ as $p\to\infty,$ where 
$Z_p=1+(\x_p^\top \Pi_p\x_p-\tr(\Pi_p))/p\geqslant 0$ a.s. Therefore, $Z_p\to 1$ in $L_1$ for any sequence of orthogonal projectors $\Pi_p,$ $p\geqslant 1$. Thus, 
\[\sup_{\Pi_p}\e|\x_p^\top \Pi_p\x_p-\tr(\Pi_p)|=o(p),\]
where the supremum is taken over all $p\times 
p $ orthogonal projectors $\Pi_p.$ Any $p\times p$ diagonal matrix $D $ with diagonal entries $\lambda_1\geqslant \ldots\geqslant \lambda_p\geqslant 0$ ($\lambda_1>0$) can be written as  \[\frac{D}{\lambda_1}=\sum_{k=1}^pw_k D_k,\]
where $\lambda_{p+1}=0,$  $w_k=(\lambda_k-\lambda_{k+1})/\lambda_1\geqslant 0$ are such that $\sum_{k=1}^pw_k=1$, and each $D_k$ is a diagonal matrix with diagonal entries in $1,\ldots,1,0,\ldots,0$ ($k$ ones). Hence, any  symmetric positive semi-definite matrix $A_p$ with $\|A_p\|=1$  can be written as a convex combination of some orthogonal projectors. As a result, by the convexity of the $L_1$-norm,
\[\sup_{ A_p}\e|\x_p^\top A_p\x_p-\tr( A_p)|\leqslant \sup_{ \Pi_p}\e|\x_p^\top \Pi_p\x_p-\tr(\Pi_p)| =o(p),\]
where $A_p$ as above. We conclude that $(\mathrm A1)$ holds. Q.e.d.

\noindent{\bf Proof of Proposition \ref{p1}.}  Suppose $(\mathrm A3)$ holds. Let us show that $(\mathrm A2)$ holds. The latter will follow from the inequality
\begin{equation}\label{inequality}
I:=\p(|\x_p^\top A_p \x_p-\tr(\Sigma_pA_p)|>\varepsilon p)\leqslant
 \frac{2\|A_p\|^2\tr(\Sigma_p^2)}{(\varepsilon p)^{2}}
\end{equation}
valid for any $\varepsilon>0$ and any $p\times p$ symmetric positive semi-definite matrix $A_p$. By Chebyshev's inequality,  $I\leqslant \var(\x_p^\top A_p \x_p)(\varepsilon p)^{-2},$ since $\e(\x_p^\top A_p \x_p)=\tr(\Sigma_pA_p)$. As a result,  we need to verify that 
\begin{equation}
\label{v1}
\var(\x_p^\top A_{p} \x_p)  \leqslant 2\|A_p\|^2\tr(\Sigma_p^2). \end{equation}

 We have $\x_p^\top A_p \x_p=\z_p^\top D_p\z_p,$ where $\z_p$ is a standard normal vector and $D_p$ is a diagonal matrix whose diagonal entries $\{\lambda_k\}_{k=1}^p$ are eigenvalues of $B_p=\Sigma_p^{1/2}A_p\Sigma_p^{1/2}$. By a direct calculation,  $\var(\z_p^\top D_{p} \z_p)=2\,\tr(D_p^2)=2\,\tr(B_p^2)$. By  Fact 1 in Appendix C and the identity $\tr( AB)=\tr(BA)$, $\tr(B_p^2)=\tr( A_p\Sigma_p A_p\Sigma_p )\leqslant \|A_p\|^2\tr(\Sigma_p^2)$. Thus, we get \eqref{inequality}.

 Assume now that $(\mathrm A2)$ holds. Let us show that $(\mathrm A3)$ holds.  Take $A_p=I_p$. Hence,
\[\frac{\x_p^\top \x_p-\tr(\Sigma_p)}{p}=\frac{\z_p^\top D_p \z_p-\tr(D_p)}{p}\pto 0,\quad p\to\infty,\]
where $\z_p$ is as above and $D_p$ is a diagonal matrix whose diagonal entries  $\{\lambda_{k}\}_{k=1}^p$ are  $\Sigma_p$'s eigenvalues arranged in descending order, i.e.  $\|\Sigma_p\|=\lambda_{1}\geqslant\ldots\geqslant \lambda_{p}\geqslant0$ and $\lambda_k=\lambda_k(p)$, $1\leqslant k\leqslant p$.

Let $\y_p=(Y_1,\ldots,Y_p)$ be an independent copy of $\z_p=(Z_1,\ldots,Z_p)$. Therefore,
$(\z_p^\top D_p \z_p-\y_p^\top D_p \y_p)/p\pto 0$
 and
\begin{align*}\e\exp\{i(\z_p^\top D_p \z_p- \y_p^\top D_p \y_p)/p\}=&\prod_{k=1}^p\e\exp\{i \lambda_{k}(Z_k^2-Y_k^2)/p\}=\prod_{k=1}^p|\varphi(\lambda_{k}/p)|^2\to 1
\end{align*}
as $p\to\infty$, where  $\varphi(t)=\e\exp\{i tZ_1^2\}$, $t\in\bR$. Hence, \[|\varphi(\|\Sigma_p\|/p)|^2=\frac{1}{|1-2i\|\Sigma_p\|/p|}=
\frac{1}{\sqrt{1+4\|\Sigma_p\|^2/p^2}}\to 1\] and $\|\Sigma_p\|/p\to0 $. As a result, 
\begin{align*}\prod_{k=1}^p|\varphi(\lambda_{k}/p)|^2
=&\prod_{k=1}^p\frac1{\sqrt{1+4\lambda_{k}^2/p^2}}
=\exp\Big\{(-2+\varepsilon_p)\sum_{k=1}^p\lambda_{k}^2/p^2\Big\}\to 1\end{align*} for some $\varepsilon_p=o(1)$.
Thus, $\tr(\Sigma_p^2)/p^2\to0,$ i.e. $(\mathrm A3)$ holds.  Q.e.d.

\noindent{\bf Proof of Proposition \ref{p3}.} For each $p\geqslant 1$, let $A_p$ be a positive semi-definite symmetric $p\times p$ matrix with $\|A_p\|=O(1)$ as $p\to\infty$. We need to show that
\begin{equation}\label{e11}
\frac{1}{p}(\x_p^\top A_p\x_p-\tr(\Sigma_pA_p))\pto0,
\end{equation}
Let further $\y_p=(Y_{p1},\ldots,Y_{pp})$ and $\z_p=(Z_{p1},\ldots,Z_{pp})$ be   centred Gaussian random vectors in $\bR^p$ with variances $A_p$ and $\Sigma_p$, respectively. Suppose also $\x_p,$ $\y_p,$ and $\z_p$ are mutually independent.

We have $|\e X_{pk}X_{pl}|\leqslant |\e X_{pk}\e[X_{pl}|\cF_k^p]|\leqslant \sqrt{\Gamma_0\Gamma_{l-k}}$ for $k\leqslant l$ and
 \[\frac{\tr(\Sigma_p^2)}{p^2}=\frac{1}{p^2}\sum_{k,l=1}^p|\e X_{pk}X_{pl}|^2\leqslant
 \frac{1}{p^2}\sum_{k,l=1}^{p}\sqrt{\Gamma_0\Gamma_{|l-k|}}\leqslant\frac2p\sum_{j=0}^p\sqrt{\Gamma_0\Gamma_j}=o(1)\]
 (recall that $\Gamma_j\to 0,$ $j\to\infty$). Thus, $(\mathrm A3)$ holds and, by Proposition \ref{p1}, \begin{equation}\label{e10}
 \frac{1}{p}(\z_p^\top A_p\z_p-\tr(\Sigma_pA_p))\pto0 .
 \end{equation} 
We need a technical lemma proved in Appendix B.
\begin{lemma}\label{l5} For each $n\geqslant 1,$ let $Z_n$ be a random variable such that $Z_n\geqslant 0$ a.s. If  $\e Z_n=O(1)$ as $n\to\infty,$ then $Z_n-\e Z_n\pto0$ iff  $\e \exp\{-Z_n\}-\exp\{-\e Z_n\}\to0 $.  
\end{lemma}

 Note that $\tr(\Sigma_p)/p\leqslant \Gamma_0$ and, by Fact 1 in Appendix C, $\tr(\Sigma_pA_p)/p=O(1)$. Thus, by $\tr(\Sigma_pA_p)=\e(\x_p^\top A_p\x_p)=\e(\z_p^\top A_p\z_p)$, Lemma \ref{l5}, and \eqref{e10}, we can prove \eqref{e11} by showing that
 \begin{align*}
 \Delta_p=\e\exp\{-\x_p^\top A_p\x_p/(2p)\}- \e\exp\{-\z_p^\top A_p\z_p/(2p)\}\to0. 
 \end{align*}
 It follows from the independence $\y_p$ and $(\x_p,\z_p)$  that  
 \[\Delta_p=\e\exp\{i(\x_p^\top\y_p)/\sqrt{p}\}- \e\exp\{i(\z_p^\top\y_p)/\sqrt{p}\},\]
 where we have used that  $\e((\x_p^\top\y_p)^2|\x_p)=\x_p^\top A_p\x_p$ and  $\e((\z_p^\top\y_p)^2|\z_p)=\z_p^\top A_p\z_p$ a.s.

Now, we will proceed in the same way as in the proof of Theorem 5 in \cite{PM}. Fix some $q,j\in \bN$ and assume w.l.o.g. that $m=p/(q+j)$ is integer (we can always add no more than $q+j$ zeros to $\x_p,\y_p,\z_p$).  Let
\[\wt\x_p=(\wt X_{1},{\bf 0}_j,\wt X_{2},{\bf 0}_j,\ldots,\wt X_{m},{\bf 0}_j)\]
for the  null vector ${\bf 0}_j$ in $\bR^j,$ where entries of $\wt X_{r}$ ($1\leqslant r\leqslant m$) are $X_{pl}$ for \[l=l_{r-1}+1,\ldots,l_{r-1} +q\quad\text{and}\quad l_{r-1}=(r-1)(q+j).\] Put also $\Delta \x_p=({\bf 0}_q,\Delta  X_{1},\ldots,{\bf 0}_q,\Delta  X_{m})$, where entries of $\Delta  X_{r}$  are $X_{pl} $ for \[l=l_{r-1}+q+1,\ldots,l_r.\] 
Define $\wt\y_p,$ $\wt\z_p$, $\Delta\y_p$, and $\Delta\z_p$ (with $\wt  Y_r,\wt  Z_r,\Delta Y_r,\Delta Z_r$) similarly. 

Since $\x_p^\top\y_p=\wt \x_p^\top\wt \y_p+(\Delta \x_p)^\top \Delta \y_p$ and  $|\exp\{ia\}-\exp\{ib\}|\leqslant |b-a|$ for $a,b\in\bR,$ we have 
\begin{align*}
|\e \exp\{ i(\x_p^\top\y_p)/\sqrt{p} \}- \e \exp\{i(\wt \x_p^\top\wt \y_p)/\sqrt{p} \} |\leqslant 
\frac{\e|(\Delta \x_p)^\top \Delta \y_p|}{\sqrt{p}}\leqslant  \frac{\sqrt{\e|(\Delta \x_p)^\top \Delta \y_p|^2}}{\sqrt{p}} 
 \end{align*}
and, by the independence of $\Delta \x_p$ and $\Delta \y_p$, 
\begin{align}\label{e16}
 \e|(\Delta \x_p)^\top \Delta \y_p|^2\leqslant  \e ( (\Delta \x_p)^\top V_p (\Delta \x_p ) )\leqslant \|V_p\| \e  (\Delta \x_p)^\top  (\Delta \x_p) \leqslant   \|A_p\| mj,
 \end{align}
 where $V_p=\var(\Delta\y_p)$ and, obviously, $\|V_p\|\leqslant \|\var(\y_p)\|=\|A_p\|$.
We can bound  $\e \exp\{i(\z_p^\top\y_p)/\sqrt{p} \}- \e \exp\{i(\wt \z_p^\top\wt \y_p)/\sqrt{p}\}$ in the same way. Combining these estimates with $m/p\leqslant 1/q$, we arrive at 
\[\Delta_p=  \e \exp\{i(\wt \x_p^\top\wt \y_p)/\sqrt{p}\}- \e \exp\{i(\wt \z_p^\top\wt \y_p)/\sqrt{p}\}+O(1)\sqrt{j/q}.\]

Fix some $\varepsilon>0$ and  set  $\wh \x_p=(\wh X_1,{\bf 0}_j,\ldots,\wh X_{m},{\bf 0}_j)$ for $\wh X_r$ having entries \[\wh X_{pl}=X_{pl}I(|X_{pl}|\leqslant \varepsilon\sqrt{p})-\e_{l_{r-1}-j}(  X_{pl}I(|X_{pl}|\leqslant \varepsilon\sqrt{p})) ,\quad l=l_{r-1}+1,\ldots,l_{r-1}+q,\]
hereinafter  $\e_l=\e\big(\cdot|\cF_{\max\{l,0\}}^p\big),$ $l\leqslant p$. 
Analogously, let $\wh \z_p=(\wh Z_1,{\bf 0}_j,\ldots,\wh Z_{m},{\bf 0}_j)$, where entries of $\wh Z_r$ are \[\wh Z_{pl}=Z_{pl}-\e^*_{l_{r-1}-j} Z_{pl},\quad l=l_{r-1}+1,\ldots,l_{r-1}+q.\]
Here $\e_l^*=\e_l^*(\cdot|\cG_{\max\{l,0\}}^p)$ for $\cG^p_{l}=\sigma(Z_{pk},1\leqslant k\leqslant l),$ $l\geqslant 1,$ and  the trivial $\sigma$-algebra $\cG_0^p$. Obviously, $\e_l^*Z_{pk}$ is a linear function in $ Z_{ps},$ $s\leqslant l$ (because $\z_p$ is a Gaussian vector) and, as a result, $\wh\z_p$ is a Gaussian vector. Also, 
\[\e \wh Z_{ps}\wh Z_{pt}=\e(\e^*_{l_{r-1}-j} \wh Z_{ps}\wh Z_{pt})= \e Z_{ps} Z_{pt}-\e(\e^*_{l_{r-1}-j} Z_{ps})(\e^*_{l_{r-1}-j} Z_{pt})\]
for $s,t=l_{r-1}+1,\ldots,l_{r-1}+q, $ and, as a result, 
\begin{equation}
\label{e13}
\var(\wt Z_r)-\var(\wh Z_r)=\e E_rE_r^\top,\quad r=1,\ldots,m, 
\end{equation}
with $E_r$ is a vector with entries $\e_{l_{r-1}-j}^*Z_{pl},$ $l=l_{r-1}+1,\ldots,l_{r-1}+q.$

We have 
\begin{align*}
|\e &\exp\{i(\wt \x_p^\top\wt \y_p)/\sqrt{p}\}- \e \exp\{i(\wt \z_p^\top\wt \y_p)/\sqrt{p}\}|\leqslant\\
&\leqslant  |\e  \exp\{i(\wt \x_p^\top\wt \y_p)/\sqrt{p}\}- \e \exp\{i(\wh \x_p^\top\wt \y_p)/\sqrt{p}\}|+\\
&\quad +|\e  \exp\{i(\wt \z_p^\top\wt \y_p)/\sqrt{p}\}- \e \exp\{i(\wh \z_p^\top\wt \y_p)/\sqrt{p}\}|\\
&\quad +|\e  \exp\{i(\wh \x_p^\top\wt \y_p)/\sqrt{p}\}- \e \exp\{i(\wh \z_p^\top\wt \y_p)/\sqrt{p}\}|. 
\end{align*}
Let us estimate the first term in the right-hand side of the last inequality. Arguing as in \eqref{e16}, we infer that, for $U_p=\var(\wt \y_p)$ (having $\|U_p\|\leqslant \|\var(\y_p)\|=\|A_p\|$), \begin{align*}
|\e &\exp\{i(\wt \x_p^\top\wt \y_p)/\sqrt{p}\}- \e \exp\{i(\wh \x_p^\top\wt \y_p)/\sqrt{p}\}|\leqslant \frac{1}{\sqrt{p}}\e|(\wt \x_p-\wh\x_p)^\top\wt \y_p|\leqslant \\
& \leqslant \frac{1}{\sqrt{p}}\big(\e (\wt \x_p-\wh\x_p)^\top U_p (\wt \x_p-\wh\x_p)\big)^{1/2}\leqslant \frac{\sqrt{\|A_p\|}}{\sqrt{p}} (\e \|\wt \x_p-\wh\x_p\|^2)^{1/2}.   
\end{align*}
Additionally, by $(a+b+c)^2\leqslant 3(a^2+b^2+c^2),$ $a,b,c\in\bR,$ 
\begin{align*}
\frac{1}{p}\e \|\wt \x_p-\wh\x_p\|^2\leqslant  3(\delta_1+\delta_2+\delta_3),
\end{align*}
 where, by \eqref{lin}, $\e X_{pl}=0 $ for all $(p,l),$
\[X_{pl}-\wh X_{pl}=X_{pl}I(|X_{pl}|> \varepsilon\sqrt{p})+\e_{l_{r-1}-j}X_{pl}-\e_{l_{r-1}-j}(  X_{pl}I(|X_{pl}|>\varepsilon\sqrt{p})),\] 
  and Jensen's inequality, we have
 \begin{align*}
&\delta_1= p^{-1}\sum \e X_{pl}^2I(|X_{pl}|> \varepsilon\sqrt{p}) \leqslant L_p(\varepsilon) =o(1),\\
&\delta_2= p^{-1}\sum\e\big|\e_{l_{r-1}-j} (X_{pl}I(|X_{pl}|> \varepsilon\sqrt{p}) )\big|^2 \leqslant \delta_1=o(1),\\
&\delta_3=  p^{-1}\sum\e\big|\e_{l_{r-1}-j} X_{pl} \big|^2 \leqslant \Gamma_j 
\end{align*}
and \begin{equation}
\label{e15}
\text{the sum }\sum=\sum_{r=1}^m\sum_{l=l_{r-1}+1}^{l_{r-1}+q}\text{  has no more than $p$ terms.}
\end{equation}
Hence, 
\begin{align} \label{e17}
\frac{1}{p}\e \|\wt \x_p-\wh\x_p\|^2\leqslant  3 \Gamma_j +o(1).
\end{align}
and $\e \exp\{i(\wt \x_p^\top\wt \y_p)/\sqrt{p}\}- \e \exp\{i(\wh \x_p^\top\wt \y_p)/\sqrt{p}\}=O(1)\sqrt{\Gamma_j}+o(1)$.

Using similar arguments, we see that
\begin{align*}
|\e  \exp\{i(\wt \z_p^\top\wt \y_p)/\sqrt{p}\}- \e \exp\{i(\wh \z_p^\top\wt \y_p)/\sqrt{p}\}|^2&\leqslant \frac{\|A_p\|}{p}\e \|\wt \z_p-\wh\z_p\|^2\\
&\quad =\frac{\|A_p\|}{p}\sum\e\big|\e_{l_{r-1}-j}^* Z_{pl}\big|^2
\end{align*}
with $\sum$ as in \eqref{e15}. Since $\x_p$ and $\z_p$ have the same covariance structure and $\z_p$ is Gaussian, then \begin{equation}
\label{e14}
\e|\e_{l_{r-1}-j}^* Z_{pl}  \big|^2\leqslant \e\big|\e_{l_{r-1}-j} X_{pl}\big|^2\leqslant \Gamma_j
\end{equation} (see Lemma 14 in \cite{PM}). Thus, 
$p^{-1}\sum\e\big|\e^*_{l_{r-1}-j} Z_{pl} \big|^2\leqslant \Gamma_j$ for
  $\sum$ from \eqref{e15}.

Combining these estimates, we deduce that 
\[|\Delta_p|\leqslant  |\e \exp\{i(\wh \x_p^\top\wt \y_p)/\sqrt{p}\}- \e \exp\{i(\wh \z_p^\top\wt \y_p)/\sqrt{p}\}|+O(1)(\sqrt{j/q}+\sqrt{\Gamma_j}).\]
Now,
\[|\e \exp\{i(\wh \x_p^\top\wt \y_p)/\sqrt{p}\}- \e \exp\{i(\wh \z_p^\top\wt \y_p)/\sqrt{p}\}|\leqslant \]
\[\leqslant\sum_{r=1}^m|\e \exp\{i(C_r+\wh X_r^\top \wt Y_r/\sqrt{p})\}- \e \exp\{i(C_r+\wh Z_r^\top \wt Y_r/\sqrt{p})\}|,\]
where, for $ r=1,\ldots, m$, \[C_r=p^{-1/2}\sum_{j<r}  \wh X_j^\top \wt Y_j+p^{-1/2}\sum_{j>r}  \wh Z_j^\top \wt Y_j\]
and the sum over the empty set is zero.

Expanding $\exp\{i(C_r+x)\}$ in Taylor's series around $x=0$, we derive 
\[\exp\{i(C_r+x)\}=\exp\{i C_r\}(1+ix-x^2/2+\theta(x))\quad\text{for}\quad|\theta(x)|\leqslant |x|^3/6\]
and $|\e \exp\{i(C_r+\wh X_r^\top \wt Y_r/\sqrt{p})\}- \e \exp\{i(C_r+\wh Z_r^\top \wt Y_r/\sqrt{p})\}|\leqslant |R_{r1}|+|R_{r2}|+|R_{r3}|,$
where
\[R_{r1}=\frac{1}{\sqrt{p}}\e \exp\{i C_r\}( \wh X_r^\top \wt Y_r-\wh Z_r^\top \wt Y_r), \]
\[R_{r2}=\frac{1}{ p}\e \exp\{i C_r\}( (\wh X_r^\top \wt Y_r)^2-(\wh Z_r^\top \wt Y_r)^2), \]
\[R_{r3}=\frac{1}{ p\sqrt{p}}\e (| \wh X_r^\top \wt Y_r|^3+|\wh Z_r^\top \wt Y_r|^3).\]
Putting $R_k=\sum_{r=1}^m R_{rk}$, $k=1,2,3,$ we see that 
\begin{equation}\label{e20}
|\Delta_p|\leqslant R_1+R_2+R_3+O(1)(\sqrt{j/q}+\sqrt{\Gamma_j}). 
\end{equation}

First, we will show that $R_1=0.$ Note that, by the definition of $\wh \x_p$,  $\e \wh X_1={\bf 0}_q $ and, for $r>1$, $\e[\wh X_r|\wh X_{r-1},\ldots,\wh  X_1]={\bf 0}_q$
a.s. In addition, by the definition of $\wh \z_p,$ $\wh Z_{r }$ and $\wh Z_{s }$ are uncorrelated when $r \neq s.$ Thus, $\wh Z_1,\ldots,\wh Z_m$ are mutually independent Gaussian vectors with mean zero.  Since $ \wh \x_p,\wt\y_p,\wh \z_p $ are also mutually independent, then, with probability one,
\[\e[\wh X_r|\wt\y_p,(\wh X_{l})_{l< r},(\wh Z_{l})_{l> r}]=\e[\wh X_r|(\wh X_{l})_{l< r}]={\bf 0}_q, \]
\[\e[\wh Z_r|\wt\y_p,(\wh X_{l})_{l< r}, (\wh Z_{l})_{l> r}]=\e \wh Z_r ={\bf 0}_q ,\]
where $r=1,\ldots,m$ and $\wh Z_{m+1}=\wh X_0={\bf 0}_q$. As a result, $R_1= 0.$ 

Let us prove that $R_3=o(1)+O(1)\varepsilon\sqrt{q}.$ Since $\wh X_r^\top \wt Y_r|\wh X_r\sim\mathcal N(0,\wh X_r^\top \wt V_r\wh X_r) $ for $\wt V_r=\var(\wt Y_r)$ (obviously, $\|\wt V_r\|\leqslant\|\var(\y_p)\|=\|A_p\|$), we   have
\[\e  |\wh X_r^\top \wt Y_r|^3= C_0\e|\wh X_r^\top \wt V_r\wh X_r|^{3/2}\leqslant 
C_0\|A_p\|^{3/2}\e\|\wh X_r\|^{3 }
\]
 for an absolute constant $C_0>0$ and each $r=1,\ldots,m$. Using the fact that entries of $\wh X_r$ are bounded by $2\varepsilon\sqrt{p}$ and the estimate  
\begin{equation}
\label{e19} 
\e\|\wh X_r\|^2\leqslant \sum_{l=l_{r-1}+1}^{l_{r-1}+q}\e X_{pl}^2I(|X_{pl}|\leqslant \varepsilon\sqrt{p}) \leqslant q\Gamma_0,
\end{equation} 
and recalling that $mq\leqslant p$, we obtain 
 \[\frac{1}{p\sqrt{p}}\sum_{r=1}^m\e  |\wh X_r^\top \wt Y_r|^3\leqslant C_0\|A_p\|^{3/2}\frac{m (q\Gamma_0) ( 4\varepsilon^2 pq )^{1/2}}{p\sqrt{p}}\leqslant 2C_0\Gamma_0\|A_p\|^{3/2} \varepsilon\sqrt{q }.\]
 
 In addition, as $\wh Z_r$ is a centred Gaussian vector whose entries have variance not greater than $\Gamma_0$ and $\big(\sum_{i=1}^q|a_i|/q\big)^{3/2}\leqslant \sum_{i=1}^q|a_i|^{3/2}/q $ for any  $a_i\in\bR$,
\[\e\|\wh Z_r\|^{3}=\e\Big|\sum_{l=l_{r-1}+1}^{l_{r-1}+q}\wh Z_{pl}^2\Big|^{3/2}\leqslant   q^{1/2}\sum_{l=l_{r-1}+1}^{l_{r-1}+q}\e|\wh Z_{pl}|^{3 }=C_0q^{3/2}\Gamma_0^{3/2}
\] for $C_0>0 $  as above. Hence, arguing as above, we infer 
\[\frac{1}{p\sqrt{p}}\sum_{r=1}^m\e  |\wh Z_r^\top \wt Y_r|^3\leqslant C_0\|A_p\|^{3/2}\frac{m(C_0q^{3/2}\Gamma_0^{3/2})}{p\sqrt{p}}\leqslant C_0^2\Gamma_0^{3/2}\|A_p\|^{3/2} \sqrt{q/p}=o(1
).\]
This proves that $R_3=o(1)+O(1)\varepsilon\sqrt{q}.$ 

To finish the proof, we need a good bound on $R_2=\sum_{r=1}^mR_{r2}$. First, note that, by the independence of  $\wh Z_r$ from everything else,
\[R_{r2}=\frac{1}{p}\e \exp\{i C_r\}( (\wh X_r^\top \wt Y_r)^2- \wt Y_r^\top \var(\wh Z_r)\wt Y_r). \]
In addition, by \eqref{e13}, Fact 1 in Appendix C, and \eqref{e14}, 
\[|\e \exp\{i C_r\}( \wt Y_r^\top \var(\wt Z_r)\wt Y_r- \wt Y_r^\top \var(\wh Z_r)\wt Y_r)|\leqslant \e ( \wt Y_r^\top (\e E_rE_r^\top ) \wt Y_r)\leqslant\]
\[\leqslant \tr( \var(\wt Y_r) \e E_rE_r^\top)\leqslant \|\var(\wt Y_r)\|\tr(\e E_rE_r^\top)\leqslant  \|A_p\|\e\|E_r\|^2=\|A_p\| q\Gamma_j.\]
Hence,  recalling  that $mq\leqslant p $  and $\var(\wt Z_r)=\var(\wt X_r)$, we get
\[R_2=\frac{1}{p}\sum_{r=1}^m \e \exp\{i C_r\}( (\wh X_r^\top \wt Y_r)^2- \wt Y_r^\top \var(\wt X_r)\wt Y_r)+O(1)\Gamma_j.\]

 Fix some integer $a>1$ and let
 \[C_r^{-a}=p^{-1/2}\sum_{j\leqslant r-a}  \wh X_j^\top \wt Y_j+p^{-1/2}\sum_{j>r}  \wh Z_j^\top \wt Y_j,\quad r=1,\ldots,m,\] where the sum over the empty set is zero. Since all entries of $\wh X_j$, $j\geqslant 1,$ are bounded by $2\varepsilon \sqrt{p}$, we have  
\begin{align} \label{e21}
 |\e (\exp\{i C_r\}&-\exp\{iC_r^{-a}\})( (\wh X_r^\top \wt Y_r)^2- \wt Y_r^\top \var(\wt X_r)\wt Y_r)|\leqslant \nonumber\\
 &\leqslant \frac{1}{\sqrt{p}} \sum_{r-a<j<r}\e |\wh X_j^\top \wt Y_j|\, |(\wh X_r^\top \wt Y_r)^2- \wt Y_r^\top \var(\wt X_r)\wt Y_r)|\nonumber\\
 &\leqslant \frac{2\varepsilon \sqrt{p}}{\sqrt{p}} \sum_{r-a<j<r}\e \|  \wt Y_j\|_1  |(\wh X_r^\top \wt Y_r)^2- \wt Y_r^\top \var(\wt X_r)\wt Y_r)|\nonumber\\
 &\leqslant 
  2\varepsilon aq^3  \sup (\e|\wh X_{ps} \wh X_{pt}|+\e|X_{ps}X_{pt}|)\e|Y_{pk}Y_{ps}Y_{pt}|\nonumber\\
  &\leqslant 
  4\varepsilon aq^3 \Gamma_0 \sup  \e|Y_{pk}|^3 \nonumber\\
   &\leqslant  C_1\varepsilon aq^3 \Gamma_0\|A_p\|^{3/2},
\,
\end{align}
where $C_1>0$ is an absolute constant, $\|x\|_1=|x_1|+\ldots+|x_q|$ for $x=(x_1,\ldots,x_q),$ and the supremum is taken over all $1\leqslant k,s,t\leqslant p$.

In addition, we have 
\[
|\e \exp\{i C_r^{-a}\}( (\wh X_r^\top \wt Y_r)^2- (\wt X_r^\top \wt Y_r)^2)|\leqslant  \e|(\wh X_r^\top \wt Y_r)^2- (\wt X_r^\top \wt Y_r)^2)|\]
and, by the Cauchy-Schwartz inequality and the independence of $(\wh X_r,\wt X_r)$ and $\wt Y_r$,
\begin{align*} 
 \e|(\wh X_r^\top \wt Y_r)^2- (\wt X_r^\top \wt Y_r)^2)|
&\leqslant  \big(\e ((\wh X_r-\wt X_r)^\top \wt Y_r)^2\big)^{1/2}\big(\e ((\wt X_r+\wh X_r)^\top \wt Y_r)^2 \big)^{1/2}\\
&\quad=\big(\e\| \wt V_r^{1/2}(\wh X_r-\wt X_r)\|^2 \big)^{1/2}\big(\e\|\wt V_r^{1/2}(\wt X_r+\wh X_r)^\top\|^2\big)^{1/2}\\
&\leqslant \|A_p\|\big(\e\|\wh X_r-\wt X_r\|^2 \big)^{1/2}\big(2\e\|\wt X_r\|^2+2\e\|\wh X_r\|^2\big)^{1/2} \end{align*}
where $\wt V_r=\var(\wt Y_r)$ has $\|\wt V_r\|\leqslant \|\var(\y_p)\|\leqslant \|A_p\|.$  Summing over $r$ and using \eqref{e17} (as well as \eqref{e19}), we arrive at 
\begin{align*} 
\frac{1}{p}\sum_{r=1}^m&\big(\e\|\wh X_r-\wt X_r\|^2 \big)^{1/2}\big(2\e\|\wt X_r\|^2+2\e\|\wh X_r\|^2\big)^{1/2}\leqslant \\
&\leqslant  \Big(\frac{1}{p}\sum_{r=1}^m \e\|\wh X_r-\wt X_r\|^2 \Big)^{1/2}\Big( \frac{2}{p}\sum_{r=1}^m( \e\|\wt X_r\|^2+ \e\|\wh X_r\|^2)\Big)^{1/2}\\
&\quad =\Big(\frac{1}{p}\e\|\wh \x_p-\wt \x_p\|^2 \Big)^{1/2}\Big( \frac{2}{p} ( \e\|\wt \x_p\|^2+ \e\|\wh \x_p\|^2)\Big)^{1/2}\\
&\leqslant \Big(3\Gamma_j+o(1)\Big)^{1/2} (2\Gamma_0)^{1/2}.
\end{align*}

Combining bounds after \eqref{e21} (and using $mq\leqslant p$), we conclude that 
\begin{align*}
|R_2|&\leqslant \frac{1}{p}\sum_{r=1}^m |D_r|+O(1)(\Gamma_j+\varepsilon aq^2)+O(1)\sqrt{\Gamma_j+o(1)}
\end{align*}
for $D_r=\e \exp\{i C_r^{-a}\}( (\wt X_r^\top \wt Y_r)^2- \wt Y_r^\top \var(\wt X_r)\wt Y_r).$
Again, the mutual independence of  $\x_p,\wt \y_p,\wh \z_p$ implies that 
\[\e\big((\wt X_r^\top \wt Y_r)^2\big| \wt\y_p,(X_{pl})_{l\leqslant (l_{r-a}-j)^+}, (\wh Z_j)_{j>r}\big)=\wt Y_r^\top \e_{l_{r-a}-j}(\wt X_r \wt X_r^\top) \wt Y_r\quad\text{a.s.}\]
and $D_r=\e e^{i C_r^{-a}}( \wt Y_r^\top [\e_{l_{r-a}-j}(\wt X_r \wt X_r^\top)- \var(\wt X_r)]\wt Y_r)$. Hence, as \[l_{r-1}+1-(l_{r-a}-j)=j+1+(a-1)(q+j)>aj,\]
we conclude that 
\begin{align*}
|D_r|&\leqslant \e | \wt Y_r^\top (\e_{l_{r-a}-j}(\wt X_r \wt X_r^\top)- \e (\wt X_r \wt X_r^\top) )\wt Y_r)|\\
&\leqslant \sum_{s,t=l_{r-1}+1}^{l_{r-1}+q} \e | Y_{ps}Y_{pt}|\e|\e_{l_{r-a}-j} (  X_{ps}   X_{qt} ) - \e ( X_{ps}  X_{pt} ) |\\&\leqslant q^2\|A_p\| \Gamma_{aj}.
\end{align*}

The latter (with $mq\leqslant p$)  implies that 
\[|R_2|\leqslant O(1)\sqrt{ o(1)+\Gamma_j}+O(1)(\Gamma_j+\varepsilon aq^2) +q\Gamma_{aj}. \]
 It follows follows from  \eqref{e20} and obtained bounds on $R_1,R_2,R_3$  that, for a sufficiently large constant $C>0,$
 \[\varlimsup_{p\to \infty} \Delta_p\leqslant  C\big(\sqrt{j/q}+\sqrt{ \Gamma_j}+\Gamma_j+\varepsilon aq^2 +q\Gamma_{aj}+\varepsilon\sqrt{q}\big)=\Delta.\]
 Tending $\varepsilon\to 0$, we see that $\Delta\to \Delta'=  C\big(\sqrt{j/q}+\sqrt{ \Gamma_j}+\Gamma_j+q\Gamma_{aj}\big)$. Then, taking $a\to\infty$ yields $\Delta'\to\Delta''=  C\big(\sqrt{j/q}+\sqrt{ \Gamma_j}\big)$. Finally, taking $q,j\to\infty$ and $j/q\to 0$, we get 
 $\Delta''\to0$. This finishes the proof of the proposition. Q.e.d.
  
\section*{Appendix B}

\noindent{\bf Proof of Lemma \ref{l4}.} We have
\begin{align}
\label{e5}
\e \frac{Z_n}{1+Z_n}-\e\frac{\e (Z_n|Y_n)}{1+\e (Z_n|Y_n)}&=\e \frac{Z_n-\e( Z_n|Y_n)}{(1+Z_n)(1+\e (Z_n|Y_n))}\nonumber\\
&=\e \frac{Z_n-\e(Z_n|Y_n)}{(1+\e(Z_n|Y_n))^2}-\e \frac{(Z_n-\e( Z_n|Y_n))^2}{(1+Z_n)(1+\e (Z_n|Y_n))^2}\nonumber\\
&=-\e \frac{(Z_n-\e( Z_n|Y_n))^2}{(1+Z_n)(1+\e (Z_n|Y_n))^2}.
\end{align}

As a result, we see that 
\[\frac{(Z_n-\e (Z_n|Y_n))^2}{(1+Z_n)(1+\e (Z_n|Y_n))^2}\pto0.\]
Since $\e Z_n$ is bounded and $Z_n\geqslant 0$ a.s., we conclude that $Z_n$  and $\e(Z_n|Y_n)$ are bounded asymptotically in probability and  $Z_n-\e(Z_n|Y_n)\pto0.$ Q.e.d.

\noindent{\bf Proof of Lemma \ref{l3}.} We have
\begin{align*} I=\frac{z_1}{1+w_1}-\frac{z_2}{1+w_2}&=\frac{z_1(1+w_2)-z_2(1+w_1)-z_1w_1+w_1z_1}{(1+w_1)(1+w_2)}\\&=
\frac{(z_1-z_2)+z_1(w_2-w_1)+w_1(z_1-z_2)}{(1+w_1)(1+w_2)}.\end{align*}
It follows from $|z_1-z_2|\leqslant \gamma$, $|w_1-w_2|\leqslant \gamma,$ and $|z_1|/|1+w_1|\leqslant M$ that 
\[|I|\leqslant \frac{\gamma(1+|z_1|+|w_1|)}{|1+w_1||1+w_2|}\leqslant \frac{\gamma}{|1+w_2||1+w_1|}+\frac{\gamma M}{|1+w_2|}+\frac{\gamma}{|1+w_2|}\,\frac{|w_1|}{|1+w_1|}.\]
In addition, we have  $|1+w_2| \geqslant \delta,$
\[|1+w_1|=|1+w_2+(w_1-w_2)|\geqslant \delta-\gamma\geqslant \delta /2,\]
\[
\frac{|w_1|}{|1+w_1|}= \dfrac{2}{|1+w_1|} I(|w_1|\leqslant 2)+
\dfrac{|w_1|}{|w_1|-1}\,I(|w_1|>2)\leqslant\begin{cases}
4/\delta,& |w_1|\leqslant 2,\\
2,& |w_1|>2.\end{cases}\]
Finally, we conclude that  $|I|\leqslant \gamma(2/\delta^2+M/\delta +4/\min\{\delta^2,2\delta\}).$ Q.e.d.

\noindent{\bf Proof of Lemma \ref{l2}.} For any given  $\varepsilon,M>0$, set
\begin{equation}\label{IM}
I_0(\varepsilon,M)=\varlimsup_{p\to\infty}\sup_{A_p}\p(|\y_p^\top A_p \y_p-\tr(\Sigma_p A_p)|>\varepsilon p),
\end{equation} 
where the supremum is taken over all  real  symmetric $p\times p$ matrices $A_p$ with  $\|A_p\|\leqslant M.$  By this definition, there are $p_j\to\infty$ and $A_{p_j}$ with $\|A_{p_j}\|\leqslant M$  such that 
\[I_0(\varepsilon,M)=\lim_{j\to\infty}\p(|\y_{p_j}^\top A_{p_j} \y_{p_j}-\tr(\Sigma_{p_j}A_{p_j})|>\varepsilon p_j).\] 
Every real symmetric matrix $A_p$ can be written as $A_p=A_{p1}-A_{p2}$ for real symmetric positive semi-definite $p\times p$ matrices $A_{pk},$ $k=1,2$, with $\|A_{pk}\|\leqslant \|A_p\|$. Moreover,
for any $\varepsilon>0$ and $p\geqslant 1$,
\begin{align} \label{ineq}
\p(|\y_p^\top A_p \y_p-\tr(\Sigma_p A_p)|>\varepsilon p)\leqslant &\sum_{k=1}^2 \p(|\y_p^\top A_{pk} \y_p-\tr(\Sigma_p A_{pk})|>\varepsilon p/2) .\end{align}
Hence, it follows from $(\mathrm A2)$ that $I_0(\varepsilon,M)=0$ for any $\varepsilon,M>0.$

If $A_p$ is any real $p\times p$ matrix and $B_p=( A_p^\top +A_p)/2$, then $\y_p^\top A_p\y_p= \y_p^\top B_p\y_p$ and, by Fact 2 in Appendix C, $\|B_p\|\leqslant \|A_p\|.$
In addition, $\tr(\Sigma_pA_p)=\tr( A_p\Sigma_p)=\tr(( A_p\Sigma_p)^\top)=\tr(\Sigma_pA_p^\top )=\tr(\Sigma_pB_p).$
Thus, if $I_1(\varepsilon,M)$ is defined as $I_0(\varepsilon,M)$ in \eqref{IM} with the supremum taken over all  real $p\times p$ matrices $A_p$ with  $\|A_p\|\leqslant M$, then  
\begin{equation}\label{I1M}
I_1(\varepsilon,M)=I_0(\varepsilon,M)=0\quad\text{ for any $\varepsilon,M>0.$}
\end{equation}

 Define  now $I_2(\varepsilon,M)$ similarly to  $I_0(\varepsilon,M)$ in \eqref{IM} with the supremum taken over all  complex $p\times p$ matrices $A_p$ with  $\|A_p\|\leqslant M$. Every such  $A_p$ can be written as $A_p=A_{p1}+iA_{p2}$ for real $p\times p$ matrices $A_{pk},$ $k=1,2.$ By Fact 3 in Appendix C,  $
\|A_{pk}\|\leqslant \|A_p\|$, $k=1,2$. Thus,  \eqref{ineq} and \eqref{I1M} yield $I_2(\varepsilon,M)=0$ for any $\varepsilon,M>0.$ Q.e.d.

\noindent{\bf Proof of Lemma \ref{l5}.} If $Z_n-\e Z_n\pto 0$ as $n\to\infty$, then $e^{\e Z_n-Z_n}-1\to 0$. Since $\e Z_n\geqslant 0$,   $e^{-\e Z_n }\leqslant 1$ and $e^{-\e Z_n}(e^{\e Z_n-Z_n}-1)= e^{-Z_n}-e^{-\e Z_n}\to 0$. By the dominated convergence theorem, $\e (e^{-Z_n}-e^{-\e Z_n})\to 0$.

Suppose now $\e e^{-Z_n}-e^{-\e Z_n}\to 0$.  By Jensen's inequality, 
\[ e^{-\e Z_n}+o(1) = \e e^{-Z_n}\geqslant  (\e e^{-Z_n/2})^2\geqslant (e^{-\e Z_n/2})^2= e^{-\e Z_n}.\]
Note also that $\e Z_n$ is bounded. Therefore, there is $c>0$ such that $e^{-\e Z_n/2}\geqslant c+o(1)$, $n\to\infty.$
Hence, $\var(e^{-Z_n/2})\to0$ as $n\to\infty$ and $\e e^{-Z_n/2}-e^{-\e Z_n/2}\to0.$ Combing these facts yields $e^{-Z_n/2}-e^{-\e Z_n/2}=e^{-\e Z_n/2}(e^{(\e Z_n-Z_n)/2}-1)\pto0$.  This shows that $e^{(\e Z_n-Z_n)/2}\pto 1$ or $Z_n-\e Z_n\pto0$. Q.e.d.

\section*{Appendix C}
In this section we list a few useful well-known  facts from linear algebra. Write $A_1\succ A_2$ for real $p\times p$ matrices $A_1,A_2$ if $A_1-A_2$ is positive semi-definite. 

 Let $B$ and $C$ be  real $p\times p$ matrices, $A$ be a complex $p\times p$ matrix, and  $z\in \mathbb C^+$.  
 \\ 
{\bf Fact 1.} $\tr(BC)\leqslant \|B\|\tr(C)$ and $\tr((BC)^2)\leqslant \|B\|^2\tr(C^2)$ if $B,C$ are symmetric and positive semi-definite.
\\
{\bf Fact 2.} If $C=(B^\top+B)/2$, then $\|C\|\leqslant \|B\|$.
\\
{\bf Fact 3.} If $A=B+iC,$ then $\|B\|\leqslant \|A\|$ and $\|C\|\leqslant \|A\|.$
\\ 
{\bf Fact 4.}  If $C$ is symmetric, then $\|(C-zI_p)^{-1}\|\leqslant 1/\Im(z) $.
\\
{\bf Fact 5.}   If  $w\in\bR^p$ and  $C$ is symmetric, then
\[ \frac{|w^\top (C-zI_p)^{-2}  w|}{|1+w^\top (C-zI_p)^{-1} w|}\leqslant \frac{1}{\Im(z)}.\]
If, in addition, $C$ is positive definite, then 
\[0\leqslant \frac{w^\top C^{-2}  w}{1+w^\top C^{-1} w}\leqslant \|C^{-1}\| \frac{w^\top C^{-1}  w}{1+w^\top C^{-1} w}\leqslant \|C^{-1}\|.\]
{\bf Fact 6.} If $B$ and $C$ are symmetric and positive semi-definite, then
\[|1+\tr(B(C-zI_p)^{-1})|\geqslant \frac{\Im(z)}{|z|}.\]
{\bf Fact 7.} If $A$ is invertible and $w\in\mathbb C^p$ satisfy $1+w^\top A^{-1}w\neq 0$, then
\[\tr(A+ww^\top )^{-1}=\tr (A^{-1})-\frac{ w^\top A^{-2}w}{1+ w^\top A^{-1}w}\quad\text{and}\quad  w^\top (A+ww^\top)^{-1}w=\frac{ w^\top A^{-1}w}{1+ w^\top A^{-1}w}.\] 
{\bf Fact 8.} If $C$ is symmetric and positive semi-definite and $\varepsilon,v>0$,  then
\[|\tr(C-(-\varepsilon+iv)I_p)^{-1}-\tr(C+\varepsilon I_p)^{-1}|\leqslant \frac{pv}{\varepsilon^2}.\] 

\noindent{\bf Proof of Fact 1.}  Since $\|B\|I_p\succ B$,   \[\|B\|C=C^{1/2}(\|B\|I_p)C^{1/2}\succ C^{1/2}BC^{1/2}\]
and $\|B\|^2C^2\succ (C^{1/2}BC^{1/2})^2$. Thus, $\tr(BC)=\tr(C^{1/2}BC^{1/2})\leqslant \|B\|\tr(C)$ and \[\tr(BCBC)=\tr((C^{1/2}BC^{1/2})^2)\leqslant \|B\|^2\tr(C^2).\]  

\noindent{\bf Proof of Fact 2.} We have  \[
\|B\|=\sqrt{\lambda_{\max}(B^\top B)}=\sqrt{\lambda_{\max}(BB^\top)}=\|B^\top\|\quad\text{and}\quad\|C\|\leqslant\frac{\|B\|+\|B^\top\|}{2}=\|B\|.\]

\noindent{\bf Proof of Fact 3.} If $A=B+iC$, then \[
\|A\|=\sup_{y\in\mathbb C^p:\,\|y\|=1}\|Ay\|\geqslant \sup_{x\in\mathbb R^p:\,\|x\|=1}\|Ax\|\geqslant
\sup_{x\in\mathbb R^p:\,\|x\|=1}\|Bx\|=\|B\|, \]
where we have used the fact that $\|Ax\|^2=\|Bx\|^2+\|Cx\|^2,$ $x\in\bR^p,$ and, for some nonzero $x_0\in\bR^p$, 
 $B^\top Bx_0=\|B\|^2x_0$ and $\|Bx_0\|=\|B\|\|x_0\|$. Similarly, we get that $\|A\|\geqslant \|C\|$. Q.e.d.

\noindent{\bf Proof of Fact 4.} The spectral norm of $A=(C-zI_p)^{-1}$ is the square root of $\lambda_{\max} ( A^*A),$ where $A^*=\overline{A^\top}=(C-\overline{z}I_p)^{-1}$. If $z=u+iv$ for $u\in \bR$ and $v=\Im(z)>0,$ then
\[A^*A=(C-\overline{z}I_p)^{-1}(C-zI_p)^{-1}=((C-uI_p)^2+v^2 I_p)^{-1}.\]
Hence, $\lambda_{\max}(A^*A)\leqslant 1/v^2$ and $\|A\|\leqslant 1/v.$ Q.e.d.

\noindent{\bf Proof of Fact 5.} Write $C=\sum_{k=1}^p\lambda_ke_ke_k^\top$ for some $\lambda_k\in\bR$ and orthonormal vectors $e_k\in\bR^p,$ $1\leqslant k\leqslant p.$ Then the result follows from the inequalities 
\begin{align*} |1+&w^\top (C-zI_p)^{-1} w|\geqslant\Im(w^\top (C-zI_p)^{-1} w)=
\Im\Big(\sum_{k=1}^p \frac{(w^\top e_k)^2}{\lambda_k-z}\Big)=\\
&=\Im (z) \sum_{k=1}^p\frac{(w^\top e_k)^2}{|\lambda_k-z|^2}\geqslant
\Im (z) \Big|\sum_{k=1}^p\frac{(w^\top e_k)^2}{(\lambda_k-z)^2}\Big|=\Im (z)|w^\top (C-zI_p)^{-2} w|.
\end{align*}
Now, if $C$ is positive definite, then $C^{-k},$ $k=1,2,$ are also positive definite and \[w^\top C^{-2} w=(C^{-1/2}w)^\top C^{-1} (C^{-1/2}w)\leqslant \|C^{-1}\|\|C^{-1/2}w\|^2=\|C^{-1}\|(w^\top C^{-1}w).\]
The latter implies desired bounds. Q.e.d.

\noindent{\bf Proof of Fact 6.} Write $z=u+iv$ for $u\in \bR$ and $v>0$. We need to prove that
\[|z||1+\tr(B(C-zI_p)^{-1})|=|z+\tr(B(C/z-I_p)^{-1})|\geqslant v.\]
Since 
\[|z+\tr(B(C/z-I_p)^{-1})|\geqslant \Im(z+\tr(B(C/z-I_p)^{-1}))=v+\Im(\tr(B(C/z-I_p)^{-1})),\]
we only need to check that $\Im\big(\tr(B(C/z-I_p)^{-1})\big)\geqslant0.$
Let
\begin{equation} \label{B}S=\Big(\frac{u}{|z|^2}C-I_p\Big)^2+\frac{v^2}{|z|^4}C^2.\end{equation}
Such $S$ is invertible, symmetric, and positive definite, since \[S=(C/z-I_p)(C/z-I_p)^*=(C/z-I_p)^*(C/z-I_p)=\]\[=
(C/\overline{z}-I_p)(C/z-I_p)=\frac{1}{|z|^2} C^2-\frac{2u}{|z|^2} C+I_p\] and $C/z-I_p=(C-zI_p)/z$ is invertible, where $A^*=\overline{A^\top}$.
Additionally, 
\begin{align*}
(C/z-I_p)^{-1}=(C/\overline{z}-I_p)S^{-1}=\Big(\frac{u}{|z|^2}C-I_p+\frac{iv}{|z|^2} C\Big)S^{-1}.
\end{align*}
Therefore,
\[\Im\big(\tr(B(C/z-I_p)^{-1})\big)=\frac{v}{|z|^2}\tr(B CS^{-1}).\]
By the definition of $S$,  $C^{1/2}$ and $S^{-1}$ commute, $CS^{-1}=C^{1/2}S^{-1}C^{1/2}$, and
\[\tr(B CS^{-1})=\tr(B C^{1/2}S^{-1}C^{1/2})=\tr(B^{1/2} C^{1/2}S^{-1}C^{1/2}B^{1/2}).\]
As it is shown above, $S$ is symmetric and positive definite. Hence, $S^{-1}$ is symmetric positive definite and 
$QS^{-1}Q^\top$ is symmetric positive semi-definite for any $p\times p$ matrix $Q$. Taking $Q=B^{1/2} C^{1/2}$ yields
 \[\tr(B CS^{-1})=\tr(QS^{-1}Q^\top)\geqslant 0.\]
 This proves the desired bound. Q.e.d. 
 
 \noindent{\bf Proof of Fact 7.} The Sherman-Morrison formula states that 
 \[(A+ww^\top )^{-1}= A^{-1}-\frac{ A^{-1}ww^\top A^{-1}}{1+ w^\top A^{-1}w}\quad\text{when}\quad 1+ w^\top A^{-1}w\neq 0.\]
 Taking traces we get the first identity. Multiplying by $w$ and $w^\top,$ we arrive at 
\[w^\top(A+ww^\top )^{-1}w= w^\top A^{-1}w-\frac{(w^\top A^{-1}w)^2}{1+ w^\top A^{-1}w}=\frac{ w^\top A^{-1}w}{1+ w^\top A^{-1}w}.\]
Q.e.d.

\noindent{\bf Proof of Fact 8.} We have 
\[ (C-(-\varepsilon+iv)I_p)^{-1}-(C+\varepsilon I_p)^{-1}=iv(C-(-\varepsilon+iv)I_p)^{-1}(C+\varepsilon I_p)^{-1}.\]
 Arguing as in the proof of Fact 4, we see that 
\[\|(C-(-\varepsilon+iv)I_p)^{-1}(C+\varepsilon I_p)^{-1}\|\leqslant 
\|(C-(-\varepsilon+iv)I_p)^{-1}\|\,\|(C+\varepsilon I_p)^{-1}\|\leqslant \varepsilon^{-2}.\]
Combining these relations and using the inequality $|\tr(A)|\leqslant p\|A\|$, we get the desired bound. Q.e.d.

\end{document}